\documentclass{article}
\usepackage{epsfig,epsfig}
\usepackage{amsmath}
\usepackage{amssymb}
\usepackage{booktabs}
\usepackage{xcolor}
\newtheorem{thm}{Theorem}[section]

\newtheorem{lemma}[thm]{Lemma}

\newtheorem{defn}[thm]{Definition}
\newtheorem{rem}[thm]{Remark}

\numberwithin{equation}{section} \topmargin=-2.5cm \oddsidemargin=0.4cm
\begin{document}
\title{Fast algorithm based on TT-M FE system for space fractional Allen-Cahn equations with smooth and non-smooth solutions
\thanks{Corresponding author.
\newline \emph{Corresponding email addresses:} mathliuyang@aliyun.com; mathliuyang@imu.edu.cn(Y.Liu).
\newline Tel. +864714991650.
\newline ORCID: https://orcid.org/0000-0001-8218-0196
\newline \emph{Manuscript submitted to Journal~~~~~~~~~~~~~~~~ May. 2018}
}}
\date{ }
\author{Baoli Yin, Yang Liu$^{*}$, Hong Li, Siriguleng He
\\\small{\emph{School of Mathematical Sciences, Inner Mongolia University, Hohhot 010021,
China;}}
}
\date{}
 \maketitle
  {\color{black}\noindent\rule[0.5\baselineskip]{\textwidth}{0.5pt} }
\noindent \textbf{Abstract:} In this article, a fast algorithm based on time two-mesh (TT-M) finite element (FE) scheme, which aims at solving nonlinear problems quickly, is considered to numerically solve the nonlinear space fractional Allen-Cahn equations with smooth and non-smooth solutions. The implicit second-order $\theta$ scheme containing both implicit Crank-Nicolson scheme and second-order backward difference method is applied to time direction, a fast TT-M method is used to increase the speed of calculation, and the FE method is developed to approximate the spacial direction. The TT-M FE algorithm includes the following main computing steps: firstly, a nonlinear implicit second-order $\theta$ FE scheme on the time coarse mesh $\tau_c$ is solved by a nonlinear iterative method; secondly, based on the chosen initial iterative value, a linearized FE system on time fine mesh $\tau<\tau_c$ is solved, where some useful coarse numerical solutions are found by the Lagrange's interpolation formula. The analysis for both stability and a priori error estimates are made in detail. Finally, three numerical examples with smooth and non-smooth solutions are provided to illustrate the computational efficiency in solving nonlinear partial differential equations, from which it is easy to find that the computing time can be saved.
\\
 {\color{black}\noindent\rule[0.5\baselineskip]{\textwidth}{0.5pt} }
\def\REF#1{\par\hangindent\parindent\indent\llap{#1\enspace}\ignorespaces}
\newcommand{\h}{\hspace{1.cm}}
\newcommand{\hh}{\hspace{2.cm}}
\newtheorem{yl}{\hspace{0.cm}Lemma}
\newtheorem{dl}{\hspace{0.cm}Theorem}
\newtheorem{re}{\hspace{0.cm}Remark}
\renewcommand{\sec}{\section*}
\renewcommand{\l}{\langle}
\renewcommand{\r}{\rangle}
\newcommand{\be}{\begin{eqnarray}}
\newcommand{\ee}{\end{eqnarray}}
\normalsize \vskip 0.2in

\section{Introduction}
\par
Fast algorithms for fractional partial differential equations (FPDEs) have been paid much attention to recently and developed rapidly. Different fast algorithms, which cover the fast computation of time fractional derivative, fast algorithm of nonlinear problem in time, fast calculation of nonlinear problem in space, fast computation of Matrix and so forth, have different acceleration strategies and features. Jiang et al. \cite{Zhangjw2} proposed a fast method of the time Caputo fractional derivative, which can reduce the computing time resulted in by the nonlocality of fractional derivative; Liu et al. \cite{Liuydyw1}, Liu et al. \cite{Liuydyw2}, and Yin et al. \cite{Liu3} considered the fast calculation for time FPDEs based on the Xu' s two-grid FE methods \cite{Xujc1}, which can reduce the calculating time yielded by the nonlinear term; Zhao et al. \cite{Zhaom1} developed a fast Hermite FE algorithm to improve the computational efficiency of Matrix, and presented a block circulant
preconditioner; Yuste and Quintana-Murillo \cite{Yuste1} presented the fast and robust adaptive methods with finite difference scheme for the time fractional diffusion equations; Xu et al. \cite{Xuqw1}, Wu and Zhou \cite{Wusl1} considered the parareal algorithms for solving the linear time fractional ordinary or partial differential equations (FO(P)DEs), respectively; Zeng et al. \cite{Zengfh2} presented a unified stable fast time-stepping method for fractional derivative and integral operators.
Recently, Liu et al. \cite{Liu1} proposed a fast TT-M FE algorithm for time fractional water wave model, which is developed to deal with time-consuming problem of nonlinear iteration used in the standard nonlinear Galerkin FE method for nonlinear term. In addition to the mentioned algorithms, there exist many other fast techniques to FPDEs or FODEs, which we will not list in this article.
\par
Here, our work is to continue developing the new fast TT-M FE algorithm \cite{Liu1} to the nonlinear space fractional Allen-Cahn problems
\begin{equation}\label{1.1}
\frac{\partial u}{\partial t}-\epsilon^2 \mathcal{L}_\alpha u+f(u)=g(\boldsymbol{z},t) ,(\boldsymbol{z},t)\in\Omega\times J,
\end{equation}
with boundary condition
\begin{equation}\label{1.2}
u(\boldsymbol{z},t)=0,  (\boldsymbol{z},t)\in\partial{\Omega}\times \bar{J},
\end{equation}
and initial condition
\begin{equation}\label{1.3}
u(\boldsymbol{z},0)=u_0(\boldsymbol{z}), \boldsymbol{z}\in\Omega,
\end{equation}
where the coefficient $\epsilon$ is a given constant, $\alpha \in(1,2)$ is the order of the fractional derivative, $f(u)=u^3-u$ is the nonlinear item, $u_0(\boldsymbol{z})$ is a given initial function, $J=(0,T]$ is the time interval with the positive constant $T$, and $\Omega=[a,b] \times [c,d]$$(\subset \mathbb{R}^2)$ is the spatial domain. With respect to fractional operator $\mathcal{L}_\alpha$, we define as follows
\begin{equation}\begin{split}\label{1.4}
\mathcal{L}_\alpha \triangleq& \mathcal{L}^x_\alpha+\mathcal{L}^y_\alpha
\end{split}\end{equation}
with
$$~\mathcal{L}^x_\alpha u \triangleq  \frac{1}{-2\cos\frac{\pi \alpha}{2}}(_{RL}D_{a,x}^{\alpha} u + _{RL}D_{x,b}^{\alpha} u),~\text{and}~\mathcal{L}^y_\alpha u \triangleq  \frac{1}{-2\cos\frac{\pi \alpha}{2}}(_{RL}D_{c,y}^{\alpha} u + _{RL}D_{y,d}^{\alpha} u),$$
where the left and right Riemann-Liouville fractional derivatives are defined respectively as
\begin{equation}\begin{split}\label{1.5}
_{RL}D_{a,x}^{\alpha} u =& \frac{1}{\Gamma(2-\alpha)}\frac{\partial^2}{\partial x^2}\int_a^x \frac{u(s,y)\mathrm{d}s}{(x-s)^{\alpha-1}},~
_{RL}D_{x,b}^{\alpha} u = \frac{1}{\Gamma(2-\alpha)}\frac{\partial^2}{\partial x^2}\int_x^b \frac{u(s,y)\mathrm{d}s}{(s-x)^{\alpha-1}};
\\_{RL}D_{c,y}^{\alpha} u =& \frac{1}{\Gamma(2-\alpha)}\frac{\partial^2}{\partial y^2}\int_c^y \frac{u(x,s)\mathrm{d}s}{(y-s)^{\alpha-1}},~
_{RL}D_{y,d}^{\alpha} u = \frac{1}{\Gamma(2-\alpha)}\frac{\partial^2}{\partial y^2}\int_y^d \frac{u(x,s)\mathrm{d}s}{(s-y)^{\alpha-1}}.
\end{split}\end{equation}
\par
The solutions' problems for space FPDEs like the equation (\ref{1.1}) have attracted a lot of people's attention. Ervin and Roop \cite{Roop1} developed the  variational solution for spatial fractional advection dispersion equations. Deng \cite{Dengwh1}, Feng et al. \cite{Fenglb1}, Bu et al \cite{Bu}, Fan et al. \cite{Fanwp1}, Zhao et al. \cite{Zhaoym1}, Li et al. \cite{Lim1}, Yue et al. \cite{YueBu}, Zhang et al. \cite{Zhanggy1}, Zhu et al. \cite{Nieyf1}, Zheng et al. \cite{Zhaozg2}, Dehghan and Abbaszadeh \cite{Dehghan1}, Chen and Wang \cite{Chenhz1}, Jin et al. \cite{Jinbt1}, Li et al. \cite{Liuzg1}
considered finite element methods for some space or space-time FPDEs. Heydari \cite{heydari1} developed the shifted Chebyshev polynomials for space fractional biharmonic equation. Bhrawy et al. \cite{Machado1},
Zeng et al. \cite{Zengfh1}, Zayernouri and Karniadakis \cite{Karniadakis1}, Zhang et al. \cite{Zhangh1} studied spectral methods for space or space-time FPDEs. Meerschaert and Tadjeran \cite{Meerschaert1}, Khaliq et al. \cite{Khaliq1}, Chen and Deng \cite{Dengwh2}, Li \cite{Licc1}, Ding and Li \cite{Dinghf1} developed some finite difference methods for space or space-time FPDEs.
  Recently, Hou et al. \cite{Tangt1} used a Crank-Nicolson finite difference methods for space fractional Allen-Cahn equations. However, the numerical studies on nonlinear space-fractional Allen-Cahn equations are still rarely considered.
\par
Here, our aim is to develop the fast TT-M FE algorithm proposed for solving the time FPDE \cite{Liu1} to solve nonlinear space fractional Allen-Cahn equations. The time direction is approximated by second-order $\theta$ scheme \cite{Liu2} derived based on the idea of the second-order $\alpha$-schemes (See Galerkin FE method by Wang, Liu et al. \cite{Wangyj1} in 2016; finite difference schemes by Gao et al. \cite{Gaogh1} in 2015 and Sun et al. \cite{Sunh} in 2016). For reducing the computing time resulted in by the existing nonlinear term, the fast TT-M FE algorithm is used to fast solve the nonlinear problem. In this article, our major work is as follows:
\\
$\blacklozenge_1$ Fast TT-M FE algorithm combined with second-order $\theta$ scheme is used to solve the nonlinear space fractional Allen Cahn equation;
 \\
$\blacklozenge_2$ Both stability and $H^{\mu}$ errors of fast TT-M FE method under the framework of second-order $\theta$ scheme are derived in detail, and some numerical examples with smooth and non-smooth data are provided to test and verify the theories.
\\
$\blacklozenge_3$ Compared with standard nonlinear Galerkin FE method, the TT-M FE algorithm can save the CPU time greatly. Moreover, it has almost the same computing accuracy as that computed by standard nonlinear Galerkin FE method.
\par
The structure of the paper is as follows. In section 2, we provide some definitions of norms and the relations between them. In section 3, we give the numerical scheme of fast TT-M FE algorithm with second-order $\theta$ scheme. In section 4, we implement the analysis of stability for the studied scheme. In section 5, we analyze the error estimates in detail. In section 6, we give some numerical examples, the analysis of the results, and the comparison between fast TT-M FE algorithm and nonlinear Galerkin FE method. Finally, we do some simple summaries for the numerical methods.
Here, we use some constants $C$, which are free of time fine mesh $\tau$, time coarse mesh $\tau_c$ and spatial mesh $h$ and may be different in different places.
\section{Preliminaries}
In this section, we state the necessary abstract setting for the analysis of the approximation to space fractional equations, which was developed by Ervin and Roop \cite{Roop1} and Roop \cite{Roop3}. Throughout, we denote $(u,v)=(u,v)_{L^2(\Omega)} =\int_{\Omega}uv \mathrm{d}\boldsymbol{z}$, $\|u\|=\|u\|_{L^2(\Omega)}=(u,u)^{1/2}$, $\mu=\frac{\alpha}{2} \in (\frac{1}{2},1)$.
\begin{defn}\label{2.1}
(Left fractional derivative space). For $\beta>0$, define the semi-norm
\begin{equation}\begin{split}\label{2.11}
|u|_{J^{\beta}_L(\Omega)}=(\|_{RL}D_{a,x}^{\beta} u\|^2+\|_{RL}D_{c,y}^{\beta} u\|^2)^{\frac{1}{2}},
\end{split}\end{equation}
and norm
\begin{equation}\begin{split}\label{2.12}
\|u\|_{J^{\beta}_L(\Omega)}=(\|u\|^2+|u|_{J_L^{\beta}(\Omega)}^2)^{\frac{1}{2}},
\end{split}\end{equation}
and denote by $J_L^{\beta}(\Omega)$(or $J_{L,0}^{\beta}(\Omega)$) the closure of $C^{\infty}(\Omega)$(or $C_0^{\infty}(\Omega)$) with respect to $\|\cdot\|_{J_L^{\beta}(\Omega)}$.
\end{defn}

\begin{defn}\label{2.2}
(Right fractional derivative space). For $\beta>0$, define the semi-norm
\begin{equation}\begin{split}\label{2.21}
|u|_{J^{\beta}_R(\Omega)}=(\|_{RL}D_{x,b}^{\alpha} u\|^2+\|_{RL}D_{y,d}^{\alpha} u\|^2)^{\frac{1}{2}},
\end{split}\end{equation}
and norm
\begin{equation}\begin{split}\label{2.22}
\|u\|_{J^{\beta}_R(\Omega)}=(\|u\|^2+|u|_{J_R^{\beta}(\Omega)}^2)^{\frac{1}{2}},
\end{split}\end{equation}
and denote by $J_R^{\beta}(\Omega)$(or $J_{R,0}^{\beta}(\Omega)$) the closure of $C^{\infty}(\Omega)$(or $C_0^{\infty}(\Omega)$) with respect to $\|\cdot\|_{J_R^{\beta}(\Omega)}$.
\end{defn}

\begin{defn}\label{2.3}
(Symmetric fractional derivative space). For $\beta>0$, $\beta \ne n-1/2$, $n \in \mathbb{N}$, define the semi-norm
\begin{equation}\begin{split}\label{2.31}
|u|_{J^{\beta}_S(\Omega)}=(|(_{RL}D_{a,x}^{\beta} u,_{RL}D_{x,b}^{\beta} u)|+|(_{RL}D_{c,y}^{\beta} u,_{RL}D_{y,d}^{\beta} u)|)^{\frac{1}{2}},
\end{split}\end{equation}
and norm
\begin{equation}\begin{split}\label{2.32}
\|u\|_{J^{\beta}_S(\Omega)}=(\|u\|^2+|u|_{J_S^{\beta}(\Omega)}^2)^{\frac{1}{2}},
\end{split}\end{equation}
and denote by $J_S^{\beta}(\Omega)$(or $J_{S,0}^{\beta}(\Omega)$) the closure of $C^{\infty}(\Omega)$(or $C_0^{\infty}(\Omega)$) with respect to $\|\cdot\|_{J_S^{\beta}(\Omega)}$.
\end{defn}

\begin{defn}\label{2.4}
(Fractional Sobolev space, see \cite{Nezza,Roop3}). For $\beta>0$, define the semi-norm
\begin{equation}\begin{split}\label{2.41}
|u|_{H^{\beta}(\Omega)}=\big\| |\xi|^{\beta} \tilde{u}(\xi) \big\|_{L^2(\mathbb{R}^2)},
\end{split}\end{equation}
and norm
\begin{equation}\begin{split}\label{2.42}
\|u\|_{H^{\beta}(\Omega)}=(\|u\|^2+|u|_{H^{\beta}(\Omega)}^2)^{\frac{1}{2}},
\end{split}\end{equation}
and denote by $H^{\beta}(\Omega)$ (or $H_0^{\beta}(\Omega)$) the closure of $C^{\infty}(\Omega)$ (or $C_0^{\infty}(\Omega)$) with respect to $\|\cdot\|_{H^{\beta}(\Omega)}$, where $\tilde{u}$ is the Fourier tansformation of u.
\end{defn}

\par
For $\Omega \subset \mathbb{R}^2$ being a convex set, the spaces $J_{L,0}^{\beta}(\Omega), J_{R,0}^{\beta}(\Omega), J_{S,0}^{\beta}(\Omega)$ and $H_0^{\beta}(\Omega)$ have the following properties, Ref.\cite{Roop1,Roop3}.

\begin{lemma}\label{2.5}
If $\beta > 0$, $\beta \ne n-\frac{1}{2}$, $n \in \mathbb{N}$, then $J_{L,0}^{\beta}(\Omega), J_{R,0}^{\beta}(\Omega), J_{S,0}^{\beta}(\Omega)$ and $H_0^{\beta}(\Omega)$ are equivalent, with equivalent seminorms and norms.
\end{lemma}

\begin{lemma}\label{2.6}
If $u \in J^{\beta}_{L,0}(\Omega)$, $0<\gamma<\beta$, then
\begin{equation}\begin{split}\label{2.61}
\|u\| \leq C|u|_{J^{\beta}_L(\Omega)}, \qquad |u|_{J^{\gamma}_L(\Omega)} \leq |u|_{J^{\beta}_L(\Omega)}.
\end{split}\end{equation}
The similar inequalities hold for $u \in J^{\beta}_{R,0}(\Omega)$ and if $\gamma \ne n-\frac{1}{2}$, $n \in \mathbb{N}$, then
\begin{equation}\begin{split}\label{2.62}
\|u\| \leq C|u|_{H^{\beta}(\Omega)}, \qquad |u|_{H^{\gamma}(\Omega)} \leq |u|_{H^{\beta}(\Omega)}.
\end{split}\end{equation}
\end{lemma}

\begin{lemma}\label{2.7}
Let $\beta>0$, $\Omega=(a,b) \times (c,d)$, $u \in J^{\beta}_{L,0}(\Omega) \cap J^{\beta}_{R,0}(\Omega)$. Then
\begin{equation}\begin{split}\label{2.71}
(_{RL}D_{a,x}^{\beta} u, _{RL}D_{x,b}^{\beta} u)=\cos(\beta \pi)\|_{RL}D_{-\infty,x}^{\beta} \bar{u}\|^2_{L^2(\mathbb{R}^2)}=\cos(\beta \pi)\|_{RL}D_{x,\infty}^{\beta} \bar{u}\|^2_{L^2(\mathbb{R}^2)},
\\(_{RL}D_{c,y}^{\beta} u, _{RL}D_{y,d}^{\beta} u)=\cos(\beta \pi)\|_{RL}D_{-\infty,y}^{\beta} \bar{u}\|^2_{L^2(\mathbb{R}^2)}=\cos(\beta \pi)\|_{RL}D_{y,\infty}^{\beta} \bar{u}\|^2_{L^2(\mathbb{R}^2)},
\end{split}\end{equation}
where $\bar{u}$ is the extension of $u$ by zero outside $\Omega$.
\end{lemma}

\begin{lemma}\label{2.8}
For any $u \in H^{\alpha}_0(\Omega)$ and $v \in H^{\mu}_0(\Omega)$ with $\mu=\alpha/2$, we have
\begin{equation}\begin{split}\label{2.81}
(_{RL}D_{a,x}^{\alpha} u, v) = (_{RL}D_{a,x}^{\mu} u, _{RL}D_{x,b}^{\mu} v), \qquad (_{RL}D_{x,b}^{\alpha} u, v) = (_{RL}D_{x,b}^{\mu} u, _{RL}D_{a,x}^{\mu} v).
\end{split}\end{equation}
\end{lemma}

\section{Numerical scheme}
To derive a fully discrete TT-M FE scheme, we first split the time
interval $[0,T]$ into a course uniform partition with the nodes $t_n=nM\tau$ $(n=0,1,2,\cdots,N)$, which satisfy $0=t_0<t_1<t_2<\cdots<t_N=T$ with the fine time step size $\tau=T/(nM)$ for some positive integer $2 \leq M \leq \frac{1}{\tau_c} $, where $\tau_c=M\tau$ is the coarse time mesh step size. Let $\psi^n=\psi(\cdot,t_n)$. Then the time-second order $\theta$ method \cite{Liu2}  for $n \ge 2$
\begin{equation}\begin{split}\label{3.1}
\mathcal{D}_{\tau}\psi(t_{n-\theta})= \frac{(3-2\theta)\psi^{n}-(4-4\theta)\psi^{n-1}+(1-2\theta)\psi^{n-2}}{2\tau}
\end{split}\end{equation}
and for the first time level we use the Crank-Nicolson discrete scheme
\begin{equation}\begin{split}\label{3.1.1}
\partial_{\frac{1}{2}} \psi=\frac{\psi^1-\psi^0}{\tau}.
\end{split}\end{equation}
\begin{lemma}\label{3.2}
For sufficiently smooth function $\psi(t)=\psi(\cdot,t) \in C^3[0,T]$ and any $\theta\in [0,\frac12]$, by Taylor expansion, the above approximation of first-order derivative at time $t_{n-\theta}$ is of second-order convergence rate, i.e.
\begin{equation}\begin{split}\label{3.21}
\psi_t(t_{n-\theta})&=\mathcal{D}_{\tau}\psi(t_{n-\theta})+R_t^{n-\theta},~n>1,
\end{split}\end{equation}
\text{and}
\begin{equation}\begin{split}\label{3.22}
\psi_t(t_{\frac{1}{2}})&=\partial_{\frac{1}{2}} \psi +E_1,~n=1,
\end{split}\end{equation}
\text{where}
\begin{equation}\begin{split}\label{3.22}
\|R_t^{n-\theta}\| \leq C \tau^2 \max_{t\in[0,T]}\|u_{ttt}\|, \quad \|E_1\| \leq C \tau^2\max_{t\in[0,T]}\|u_{ttt}\|,
\end{split}\end{equation}
\text{with the constant $C$ independent of $\tau$}.
\end{lemma}

\par To formulate the time semidiscrete scheme and stability, we state the following lemmas with respect to time $t_{n-\theta}$.

\begin{lemma}(See \cite{Liu2})\label{3.3}
For sufficiently smooth function $\psi(t)=\psi(\cdot,t) \in C^2[0,T]$ and function $f(t) \in C^2[0,T]$, at time $t_{n-\theta}$, the following approximate formula
\begin{equation}\begin{split}\label{3.31}
\psi(t_{n-\theta})=&(1-\theta)\psi(t_{n})+\theta \psi(t_{n-1})+E_2^{n-\theta},
\\ f(\psi(t_{n-\theta}))=&(1-\theta)f(\psi(t_{n}))+\theta f(\psi(t_{n-1}))+E_3^{n-\theta},
\end{split}\end{equation}
holds for any $\theta\in [0,\frac12]$, $n \ge 1$, where $|E_2^{n-\theta}| \leq C \tau^2$ and $|E_3^{n-\theta}| \leq C \tau^2$ with constant $C$ independent of $\tau$. We take the following notations
\begin{equation}\begin{split}\label{3.32}
\psi^{n-\theta}& = (1-\theta)\psi^{n}+\theta \psi^{n-1},
\\f^{n-\theta}(\psi)& = (1-\theta)f(\psi^{n})+\theta f(\psi^{n-1}).
\end{split}\end{equation}
\end{lemma}

\begin{lemma}(See \cite{Liu2})\label{3.4}
For series $\{\psi^n\}$ and $0\leq\theta\leq1/2$, the following inequalities hold
\begin{equation}\label{3.41}\begin{split}
\Big{(}\mathcal{D}_{\tau}\psi^{n-\theta},\psi^{n-\theta}\Big{)}&\geq \frac{1}{4\tau_c}(\mathbb{H}[\psi^{n}]-\mathbb{H}[\psi^{n-1}])
, \quad n \geq 2,
\\\mathbb{H}[\psi^{n}]&\geq\frac{1}{1-\theta}\|\psi^{n}\|^{2}, \quad n \geq 2,
\end{split}\end{equation}
where $\mathbb{H}[\psi^{n}]=(3-2\theta){\|\psi^{n}\|}^{2}-(1-2\theta){\|\psi^{n-1}\|}^{2}+(2-\theta)(1-2\theta){\|\psi^{n}-\psi^{n-1}\|}^{2},~n\geq1$.
\end{lemma}
\textbf{Proof.} See the results in \cite{Liu1}. Also follow the related results based on $\alpha$-scheme in \cite{Wangyj1,Gaogh1,Sunh,Liun1} to easily get the conclusion.
\par
Using the above $\theta$ method and Lemma \ref{2.8}, the temporal semidiscrete scheme of (\ref{1.1})-(\ref{1.3}) is to find $u^n:[0,T] \longmapsto H^{\mu}_0$, the closure of $C_0^{\infty}$ with respect to $\|\cdot\|_{H^{\mu}(\Omega)}$, which satisfies for any $v \in H^{\mu}_0$
\par
Case $n=1$:
\begin{equation}\label{3.5}\begin{split}
\Big{(}\frac{u^1-u^0}{\tau},v\Big{)}+B\Big{(}\frac{u^0+u^1}{2},v\Big{)}+\Big{(}\frac{f(u^0)+f(u^1)}{2},v\Big{)}=\Big{(}\frac{g^0+g^1}{2},v\Big{)},
\end{split}
\end{equation}
\par
Case $n>1$:
\begin{equation}\label{3.6}\begin{split}
(\mathcal{D}_{\tau}u^{n-\theta},v)+B(u^{n-\theta},v)+(f^{n-\theta}(u),v)=(g^{n-\theta},v),
\end{split}
\end{equation}
with $u^0=u_0(x,y)$. Here the bilinear form $B(u,v)$ is defined as
\begin{equation}\label{3.7}\begin{split}
B(u,v)=\frac{\epsilon^2}{2 \cos\pi \mu}\big( (_{RL}D_{a,x}^{\mu} u,_{RL}D_{x,b}^{\mu} v)+(_{RL}D_{x,b}^{\mu} u,_{RL}D_{a,x}^{\mu} v) \\+(_{RL}D_{c,y}^{\mu} u,_{RL}D_{y,d}^{\mu} v)+(_{RL}D_{y,d}^{\mu} u,_{RL}D_{c,y}^{\mu} v)\big).
\end{split}
\end{equation}
\par
By Lemmas \ref{2.5}-\ref{2.7}, we have the crucial properties of the bilinear form $B(u,v)$.

\begin{thm}\label{3.8}(See \cite{Bu})
The bilinear form $B(u,v):H^{\mu}_0 \times H^{\mu}_0 \to \mathbb{R}$ is continuous and coercive, i.e. there is a constant C independent of $u$ and $v$ such that
\begin{equation}\label{3.81}\begin{split}
|B(u,v)| \leq C \|u\|_{H^{\mu}(\Omega)}\|v\|_{H^{\mu}(\Omega)}
\end{split}
\end{equation}
and
\begin{equation}\label{3.82}\begin{split}
B(u,u) \ge C \|u\|_{H^{\mu}(\Omega)}^2.
\end{split}
\end{equation}
\end{thm}

\par
To formulate FE scheme, we define $V_h$ as the subspace of $H_0^{\mu}$, i.e.
\begin{equation}\label{3.83}\begin{split}
V_h=\{v \in H^{\mu}_0(\Omega):v|_{e} \in P_k(x,y)\},
\end{split}
\end{equation}
where $P_k(x,y)$ is the set of linear polynomials of $x,y$ with the degree no greater than $k \in \mathbb{Z}^{+}$. Then the complete discrete scheme of (\ref{1.1})-(\ref{1.3}) is to find $U^n:[0,T] \longmapsto V_h$, such that for any $v_h \in V_h$
\par
Case $n=1$:
\begin{equation}\label{3.91}\begin{split}
\Big(\frac{U^1-U^0}{\tau},v_h\Big)+B\Big(\frac{U^0+U^1}{2},v_h\Big)+\Big(\frac{f(U^0)+f(U^1)}{2},v_h\Big)=\Big(\frac{g^0+g^1}{2},v_h\Big),
\end{split}
\end{equation}
\par
Case $n>1$:
\begin{equation}\label{3.101.}\begin{split}
(\mathcal{D}_{\tau}U^{n-\theta},v_h)+B(U^{n-\theta},v_h)+(f^{n-\theta}(U),v_h)=(g^{n-\theta},v_h),
\end{split}
\end{equation}
with $U^0=u_{h0}(x,y)$, a proper approximation of $u_0(x,y)$.
Using the formulation (\ref{3.1}) and (\ref{3.32}), the expension form of (\ref{3.101.}) is as follows
\begin{equation}\label{3.11.}\begin{split}
\frac{1}{2\tau}\big((3-2\theta)U^n-(4-4\theta)U^{n-1}+(1-2\theta)U^{n-2}, v_h \big)+(1-\theta)B(U^n,v_h)+\theta B(U^{n-1},v_h)
\\+(1-\theta)(f(U^n),v_h)+\theta (f(U^{n-1}),v_h)=(1-\theta)(g^n,v_h)+\theta (g^{n-1},v_h).
\end{split}
\end{equation}
Due to the nonelinear item $(f(U^n),v_h)$, we solve the equation by a iteration method on the time mesh $\tau$.
\par
To improve the computation efficiency of the FE discrete system (\ref{3.91}) and (\ref{3.101.}), we consider the following TT-M system based on FE method, which includes the time coarse mesh $\tau_c$ and the time fine mesh $\tau$, see Ref.\cite{Liu1}
\\
\textbf{Step I:} Firstly, we get the coarse time mesh numerical approximation $U_C^n$ by equations (\ref{3.91}) and (\ref{3.101.}), i.e.
\par
Case $n=1$:
\begin{equation}\label{3.9}\begin{split}
\Big(\frac{U_C^1-U_C^0}{\tau_c},v_h\Big)+B\Big(\frac{U_C^0+U_C^1}{2},v_h\Big)+\Big(\frac{f(U_C^0)+f(U_C^1)}{2},v_h\Big)=\Big(\frac{g^0+g^1}{2},v_h\Big),
\end{split}
\end{equation}
\par
Case $n>1$:
\begin{equation}\label{3.10.}\begin{split}
(\mathcal{D}_{\tau_c}U_C^{n-\theta},v_h)+B(U_C^{n-\theta},v_h)+(f^{n-\theta}(U_C),v_h)=(g^{n-\theta},v_h),
\end{split}
\end{equation}
with $U_C^0=u_{h0}(x,y)$, a proper approximation of $u_0(x,y)$.
\\
\textbf{Step II:} Secondly, we use Lagrange's interpolation formula to get the values between $U_C^n$ and $U_C^{n+1} (n=0,1,2,\cdots,N-1)$ based on the fine time mesh. Denote by $U_I^{m}$ the interpolated results where $m=0,1,\cdots,M,M+1,\cdots,2M,\cdots,NM$ is the fine time mesh partition.
\\
\textbf{Step III:} Finally, based on the solution $U_I^m\in V_h$, the following linear system on the fine time mesh $\tau$ is considered to find $U_F^m:[0,T]\longmapsto V_h$ such that for any $v_h \in V_h$
\par
Case $m=1$:
\begin{equation}\label{3.12.}\begin{split}
\Big(\frac{U_F^1-U_F^0}{\tau},v_h\Big)+B\Big(\frac{U_F^0+U_F^1}{2},v_h\Big)+\frac{1}{2}\Big(f(U_I^1)+(U_F^1-U_I^1)f_u(U_I^1),v_h\Big)
\\ + \frac{1}{2}\Big(f(U_F^0),v_h \Big)=\Big(\frac{g^0+g^1}{2},v_h\Big),
\end{split}
\end{equation}
\par
Case $m>1$:
\begin{equation}\label{3.13.}\begin{split}
\Big(\mathcal{D}_{\tau}U_F^{m-\theta},v_h\Big)+B\Big(U_F^{m-\theta},v_h\Big)+(1-\theta)\Big(f(U_I^m)+(U_F^m-U_I^m)f_u(U_I^m),v_h\Big)
\\+\theta\Big(f(U_F^{m-1}),v_h\Big)=\Big(g^{m-\theta},v_h\Big),
\end{split}
\end{equation}
\\where $f_u$ is the derivative of $u$.
\begin{rem}\label{rem1}

\emph{\textbf{(i)}}. Fast TT-M algorithm is proposed by solving the time fractional PDE in \cite{Liu1}, and here for the first time applied to the nonlinear space fractional PDE.
\\
\emph{\textbf{(ii)}}. Direct computing the nonlinear system (\ref{3.91})-(\ref{3.11.}) by iteration is a time consuming work. So we consider TT-M FE algorithm to solve the nonlinear space fractional problem. We will show these comparisons of computing time in section 6.
\\
\emph{\textbf{(iii)}}.
Here, we first combine second-order $\theta$-scheme \cite{Liu2} with TT-M FE algorithm to nonlinear space fractional problem. Compared with the linearized $\theta$ scheme in \cite{Liu2}, we use the nonlinear $\theta$-scheme in time.
\end{rem}
\section{The analysis of stability}
\begin{thm}\label{4.1.}
For the coarse time mesh system (\ref{3.9})-(\ref{3.10.}), the following stable inequality holds
\begin{equation}\label{4.1.1}\begin{split}
\|U_C^n\|^2 \leq C\|U_C^0\|^2 + C \tau_c \sum_{k=0}^n \|g^k\|^2.
\end{split}
\end{equation}
\\Further, the stability of TT-M scheme (\ref{3.12.})-(\ref{3.13.}) holds
\begin{equation}\label{4.1.2}\begin{split}
\|U_F^m\|^2 \leq C(\|U_C^0\|^2+\|U_F^0\|^2) +C \tau\sum_{k=0}^{m+M} \|g^k\|^2.
\end{split}
\end{equation}
\end{thm}
\textbf{Proof.} \textbf{Step I:} For (\ref{4.1.1}), we substitute $U_C^{n-\theta}$ for $v_h$ in (\ref{3.10.}), i.e.
\begin{equation}\label{4.1.3}\begin{split}
(\mathcal{D}_{\tau_c}U_C^{n-\theta},U_C^{n-\theta})+B(U_C^{n-\theta},U_C^{n-\theta})+(f^{n-\theta}(U_C),U_C^{n-\theta})=(g^{n-\theta},U_C^{n-\theta}).
\end{split}
\end{equation}
By the coercivity of the bilinear form $B(u,v)$, we have the following inequality
\begin{equation}\label{4.1.4}\begin{split}
(\mathcal{D}_{\tau_c}U_C^{n-\theta},U_C^{n-\theta})+(f^{n-\theta}(U_C),U_C^{n-\theta}) \leq (g^{n-\theta},U_C^{n-\theta}).
\end{split}
\end{equation}
For any $u,v \in L^2(\Omega)$, by Cauchy-Schwarz inequality and Young inequality, one easily gets
\begin{equation}\label{4.1.5}\begin{split}
|(u^{n-\theta},v^{n-\theta})| \leq \frac{1}{2}(1-\theta)(\|u^n\|^2+\|v^n\|^2)+\frac{\theta}{2}(\|u^{n-1}\|^2+\|v^{n-1}\|^2).
\end{split}
\end{equation}
Using Lemma \ref{3.4} as well as the above inequality, (\ref{4.1.4}) is formulated as
\begin{equation}\label{4.1.6}\begin{split}
\frac{1}{4 \tau_c}(\mathbb{H}[U_C^n]-\mathbb{H}[U_C^{n-1}]) \leq \frac{1-\theta}{2}\|g^n\|^2+\frac{\theta}{2}\|g^{n-1}\|^2+(1-\theta)\|U_C^n\|^2+\theta \|U_C^{n-1}\|^2+\mathcal F^n,
\end{split}
\end{equation}
where $\mathcal F^n = \frac{1}{2}(1-\theta)\|f(U_C^n)\|^2 + \frac{\theta}{2} \|f(U_C^{n-1})\|^2 \leq C(\|U_C^n\|^2+\|U_C^{n-1}\|^2)$ with the constant $C$ independent of $n$.
\\Add up the inequality (\ref{4.1.6}) from $2$ to $n$, then
\begin{equation}\label{4.1.7}\begin{split}
\mathbb{H}[U_C^n]-\mathbb{H}[U_C^1] &\leq  C \tau_c\sum_{k=2}^n(\|g^k\|^2+\|g^{k-1}\|^2)+C \tau_c \sum_{k=2}^n(\|U_C^k\|^2+ \|U_C^{k-1}\|^2)
\\&\leq  C \tau_c\sum_{k=1}^n \|g^k\|^2+C \tau_c \sum_{k=1}^n \|U_C^k\|^2.
\end{split}
\end{equation}
Again using Lemma \ref{3.4}, we can get
\begin{equation}\label{4.1.8}\begin{split}
\frac{1}{1-\theta}\|U_C^n\|^2 \leq \mathbb{H}[U_C^n] \leq  \mathbb{H}[U_C^1]+ C \tau_c\sum_{k=1}^n \|g^k\|^2+C \tau_c \sum_{k=1}^n \|U_C^k\|^2.
\end{split}
\end{equation}
Here, $\mathbb{H}[U_C^1]=(3-2\theta)\|U_C^1\|^2-(1-2\theta)\|U_C^0\|^2+(2-\theta)(1-2\theta)\|U_C^1-U_C^0\|^2$, by triangle inequality and Yong inequality, we get
\begin{equation}\label{4.1.9}\begin{split}
\mathbb{H}[U_C^1] \leq C\|U_C^1\|^2+C\|U_C^0\|^2.
\end{split}
\end{equation}
To estimate $\|U_C^1\|^2$, substitute $(U_C^0+U_C^1)/2$ for $v_h$ in (\ref{3.9}), i.e.
\begin{equation}\label{4.1.10.}\begin{split}
\Big(\frac{U_C^1-U_C^0}{\tau_c},\frac{U_C^0+U_C^1}{2}\Big)+B\Big(\frac{U_C^0+U_C^1}{2},\frac{U_C^0+U_C^1}{2}\Big)+\Big(\frac{f(U_C^0)+f(U_C^1)}{2},\frac{U_C^0+U_C^1}{2}\Big)
\\=\Big(\frac{g^0+g^1}{2},\frac{U_C^0+U_C^1}{2}\Big).
\end{split}
\end{equation}
Using the similar analysis we have
\begin{equation}\label{4.1.11.}\begin{split}
\|U_C^1\|^2 \leq C\|U_C^0\|^2 + C \tau_c (\|g^0\|^2+\|g^1\|^2).
\end{split}
\end{equation}
Combining (\ref{4.1.8}), (\ref{4.1.9}) and (\ref{4.1.11.}), the following inequality holds for sufficiently small $\tau_c$
\begin{equation}\label{4.1.12.}\begin{split}
\|U_C^n\|^2 \leq C\|U_C^0\|^2 + C \tau_c \sum_{k=0}^n \|g^k\|^2 + C \tau_c \sum_{k=0}^n \|U_C^k\|^2.
\end{split}
\end{equation}
Then the discrete Gronwall inequality shows (\ref{4.1.1}).
\\\textbf{Step II:} For (\ref{4.1.2}), we substitute $U_F^{m-\theta}$ for $v_h$ in (\ref{3.13.}), i.e.
\begin{equation}\label{4.1.13.}\begin{split}
(\mathcal{D}_{\tau}U_F^{m-\theta},U_F^{m-\theta})+B(U_F^{m-\theta},U_F^{m-\theta})+(1-\theta)(f(U_I^m)+(U_F^m-U_I^m)f_u(U_I^m),U_F^{m-\theta})
\\+\theta(f(U_F^{m-1}),U_F^{m-\theta})=(g^{m-\theta},U_F^{m-\theta}).
\end{split}
\end{equation}
Using the techniques applied to (\ref{4.1.1}), we easily get the inequality
\begin{equation}\label{4.1.14.}\begin{split}
\mathbb{H}[U_F^m]-\mathbb{H}[U_F^1] &\leq  C \tau\sum_{k=2}^m(\|g^k\|^2+\|g^{k-1}\|^2)+C \tau \sum_{k=2}^m(\|U_F^k\|^2+ \|U_F^{k-1}\|^2)+C \tau \sum_{k=1}^m \|U_I^k\|^2
\\&\leq  C \tau\sum_{k=1}^m \|g^k\|^2+C \tau \sum_{k=1}^m \|U_F^k\|^2+C \tau \sum_{k=1}^m \|U_I^k\|^2.
\end{split}
\end{equation}
Now we estimate the Lagrange interpolation item $\|U_I^k\|$. Denote by $n=\lceil \frac{k}{M} \rceil$, the smallest integer that is equal to or greater than $\frac{k}{M}$, then by interpolation formula we have
\begin{equation}\label{4.1.15.}\begin{split}
U^k_I=\lambda_k U^{n-1}_C+(1-\lambda_k)U^{n}_C
\end{split}
\end{equation}
where $\lambda_k=n-\frac{k}{M} \in [0,1)$.
\begin{equation}\label{4.1.16.}\begin{split}
\tau\sum_{k=1}^m \|U_I^k\|^2 &\leq  \tau\sum_{k=1}^m \| \lambda_k U^{n-1}_C+(1-\lambda_k)U^{n}_C \|^2
\\ &\leq  C\tau \sum_{k=1}^m (\| U^{n-1}_C \|^2 +\| U^{n}_C \|^2) \leq C \tau\sum_{k=1}^{M \lceil \frac{m}{M} \rceil} (\| U^{n-1}_C \|^2 +\| U^{n}_C \|^2)
\\ &\leq C \tau\sum_{l=0}^{\lceil \frac{m}{M} \rceil-1} \sum_{k=1+lM}^{(l+1)M} (\| U^{n-1}_C \|^2 +\|U^{n}_C \|^2)
\\ &=  C \tau\sum_{l=0}^{\lceil \frac{m}{M} \rceil-1} \sum_{k=1+lM}^{(l+1)M} (\| U^{l}_C \|^2 +\|U^{l+1}_C \|^2) = C M\tau\sum_{l=0}^{\lceil \frac{m}{M} \rceil-1} (\| U^{l}_C \|^2 +\|U^{l+1}_C \|^2)
\\&\leq  C \tau_c\sum_{l=0}^{n} \| U^{l}_C \|^2 \leq  C \tau_c\sum_{l=0}^{n} (C\|U_C^0\|^2 + C \tau_c \sum_{k=0}^l \|g^k\|^2)
\\ &\leq  C\|U_C^0\|^2+C\tau_c^2\sum_{l=0}^n \sum_{k=0}^l \|g^k\|^2 = C\|U_C^0\|^2+C\tau_c^2\sum_{k=0}^n \sum_{l=k}^n \|g^k\|^2
\\ &\leq C\|U_C^0\|^2+C\tau_c\sum_{k=0}^n \|g^k\|^2 \leq  C\|U_C^0\|^2+C\tau\sum_{k=0}^{m+M} \|g^k\|^2.
\end{split}
\end{equation}
With (\ref{4.1.14.}) and (\ref{4.1.16.}) the following inequality holds
\begin{equation}\label{4.1.17.}\begin{split}
\mathbb{H}[U_F^m] \leq  \mathbb{H}[U_F^1] + C\|U_C^0\|^2 +  C \tau\sum_{k=0}^{m+M} \|g^k\|^2+C \tau \sum_{k=1}^m \|U_F^k\|^2.
\end{split}
\end{equation}
To estimate $\mathbb{H}[U_F^1]$ we substitute $(U_F^0+U_F^1)/2$ for $v_h$ in (\ref{3.12.}) and use the quite similar analysis above-mentioned to complete the proof for (\ref{4.1.2}). The proof is completed.
\section{Error analysis}
\par
Firstly, we give some lemmas and definitions for later analysis. Define $|u|_{\mu}=B(u,u)^{1/2}$ and $\|u\|_{\mu}=(\|u\|^2+|u|^2_{\mu})^{1/2}$  for $u \in H_0^{\mu}$.
\begin{defn}\label{5.2}
The orthogonal projection operator $P_h : H^{\mu}_0 \to V_h$ is defined as
\begin{equation}\label{5.2.1}\begin{split}
B(u-P_h u,v)=0, \qquad u\in H^{\mu}_0(\Omega) \quad \forall v \in V_h.
\end{split}
\end{equation}
\end{defn}

\begin{lemma}\label{5.3}(See  \cite{Bu})
Let s and r be real numbers satisfying $0<r \leq k+1$, $0 \leq s < r$. Then there exist a projector $\Pi_h$ and a positive constant $C$ depending only on $\Omega$ such that, for any function $u \in H^s(\Omega)$, the following estimate holds
\begin{equation}\label{5.3.1}\begin{split}
\|u- \Pi_h u\|_{H^s(\Omega)} \leq Ch^{r-s}\|u\|_{H^r(\Omega)}.
\end{split}
\end{equation}
\end{lemma}
\par
By Lemma \ref{5.3}, the operator $P_h$ defined in (\ref{5.2.1}) has the following estimate property with respect to the seminorm $|\cdot|_{\mu}$.
\begin{lemma}\label{5.4}(See \cite{Zengfh1})
Let $\mu$ and r be  real numbers satisfying $1/2 < \mu < 1$, $\mu < r \leq k+1$. Then there exists  a positive constant $C$ independent of $h$ such that, for any function $u \in H^r(\Omega) \cap H^{\mu}_0(\Omega)$, the following estimate holds
\begin{equation}\label{5.4.1}\begin{split}
\|u- P_h u\|_{\mu} \leq Ch^{r-\mu}\|u\|_{H^r(\Omega)}.
\end{split}
\end{equation}
\end{lemma}
\par
Similar to Lemma \ref{3.4}, we have the following estimates,
\begin{lemma}\label{5.6}
For series $\{\phi^n \}$ and $0\leq\theta\leq1/2$, the following inequalities hold
\begin{equation}\label{3.41}\begin{split}
B\big{(}\mathcal{D}_{\tau}\phi^{n-\theta},\phi^{n-\theta}\big{)}&\geq   \frac{1}{4\tau}(\mathbb{L}[\phi^{n}]-\mathbb{L}[\phi^{n-1}]), \quad n \geq 2,
\\\mathbb{L}[\phi^{n}]&\geq  \frac{1}{1-\theta}|\phi^{n}|^2_{\mu},  \quad n \geq 2,
\end{split}\end{equation}
where $\mathbb{L}[\phi^{n}]=(3-2\theta){|\phi^{n}|^2_{\mu}}-(1-2\theta){|\phi^{n-1}|^2_{\mu}}+(2-\theta)(1-2\theta){|\phi^{n}-\phi^{n-1}|^2_{\mu}},~n\geq1$.
\end{lemma}

\begin{thm}\label{5.7}
Suppose $u$, $U_C^n$, $U_F^m$, are the solutions of initial problem (\ref{1.1})-(\ref{1.3}), the coarse time mesh problem (\ref{3.9})-(\ref{3.10.}) and the fine time mesh problem (\ref{3.12.})-(\ref{3.13.}), respectively, with the assumption $u \in C^3(0,T;H^{k+1}(\Omega))$.Let $\mu$ and r be real numbers satisfying $\mu<r \leq k+1$, $\frac{1}{2} < \mu < 1$. Then there exits a positive constant $C$ independent of $\tau_c$, $\tau$, and $h$ such that
\begin{equation}\label{5.7.1}\begin{split}
|u(t_n)-U_C^n|_{\mu} \leq  C(\tau_c^2 + h^{r-\mu}),
\end{split}
\end{equation}
\begin{equation}\label{5.7.2}\begin{split}
\|u(t_m)-U_F^m\|_{\mu} \leq  C(\tau_c^4 + \tau^2 + h^{r-\mu}).
\end{split}
\end{equation}
\end{thm}
\textbf{Proof.} The weak formula of the initial system (\ref{1.1}) is for any $v \in H^{\mu}_0(\Omega)$
\par
Case $n=1$:
\begin{equation}\label{5.7.3}\begin{split}
&(u_t(t_{\frac{1}{2}}),v)+B(u(t_\frac{1}{2}),v)+(f(u(t_\frac{1}{2})),v)=(g(t_\frac{1}{2}),v),~or
\\&(\partial_{\frac{1}{2}}u,v)+B(u(t_\frac{1}{2}),v)+(f^{\frac{1}{2}}(u),v)=(g^{\frac{1}{2}},v)-(E_1-E_2^{\frac{1}{2}}+E_3^{\frac{1}{2}},v),
\end{split}
\end{equation}
\par
Case $n>1$:
\begin{equation}\label{5.7.4}\begin{split}
&(u_t(t_{n-\theta}),v)+B(u(t_{n-\theta}),v)+(f(u(t_{n-\theta})),v)=(g(t_{n-\theta}),v),~ or
\\&(\mathcal{D}_{\tau_c}u^{n-\theta},v)+B(u(t_{n-\theta}),v)+(f^{n-\theta}(u),v)=(g^{n-\theta},v)-(R_t^{n-\theta}-E_2^{n-\theta}+E_3^{n-\theta},v),
\end{split}
\end{equation}
where $u_t(t_{\frac{1}{2}})=\partial_{\frac{1}{2}}u+E_1$, $f(u(t_{n-\theta}))=f^{n-\theta}(u)+E_3^{n-\theta}$, and $g(t_{n-\theta})=g^{n-\theta}+E_2^{n-\theta}$.
\\\textbf{Step I:} For (\ref{5.7.1}), let $U_C^n-u(t_n) = (U_C^n- P_h u(t_n)) + (P_h u(t_n)- u(t_n)) \triangleq \xi_c^n + \rho_c^n$. First subtracting (\ref{5.7.4}) from (\ref{3.10.}) we have
\begin{equation}\label{5.7.5}\begin{split}
(\mathcal{D}_{\tau_c}\xi_c^{n-\theta},v_h)+&B(\xi_c^{n-\theta},v_h)+(f^{n-\theta}(U_C)-f^{n-\theta}(u),v_h)
\\&=(R_t^{n-\theta}-\mathcal{D}_{\tau_c}\rho_c^{n-\theta}-E^{n-\theta}_2+E^{n-\theta}_3,v_h), \quad \forall v_h \in V_h.
\end{split}
\end{equation}
Choosing $v_h=\mathcal{D}_{\tau_c}\xi_c^{n-\theta}$ in (\ref{5.7.5}), we have the following estimate by  Cauchy-Schwarz inequality and Young inequality
\begin{equation}\label{5.7.6}\begin{split}
B(\xi_c^{n-\theta},\mathcal{D}_{\tau_c}\xi_c^{n-\theta}) \leq \frac{1}{2}\|f^{n-\theta}(U_C)-f^{n-\theta}(u)\|^2 + \frac{1}{2}\|R_t^{n-\theta}-\mathcal{D}_{\tau_c}\rho_c^{n-\theta}-E^{n-\theta}_2+E^{n-\theta}_3\|^2.
\end{split}
\end{equation}
By Lemma \ref{5.6} we have
\begin{equation}\label{5.7.7}\begin{split}
\mathbb{L}[\xi_c^{n}]-\mathbb{L}[\xi_c^{n-1}] &\leq 8\tau_c(\|R_t^{n-\theta}\|^2+\|\mathcal{D}_{\tau_c}\rho_c^{n-\theta}\|^2+\|E_2^{n-\theta}\|^2+\|E_3^{n-\theta}\|^2)+2\tau_c\|f^{n-\theta}(U_C)-f^{n-\theta}(u)\|^2
\\& \leq 8\tau_c(\|R_t^{n-\theta}\|^2+\|\mathcal{D}_{\tau_c}\rho_c^{n-\theta}\|^2+\|E_2^{n-\theta}\|^2+\|E_3^{n-\theta}\|^2)
\\&\quad+4(1-\theta)^2\tau_c\|f(U^n_C)-f(u^n)\|^2+4\theta^2\tau_c\|f(U^{n-1}_C)-f(u^{n-1})\|^2
\\& \leq 8\tau_c(\|R_t^{n-\theta}\|^2+\|\mathcal{D}_{\tau_c}\rho_c^{n-\theta}\|^2+\|E_2^{n-\theta}\|^2+\|E_3^{n-\theta}\|^2)
\\&\quad+C\tau_c\|U^n_C-u(t_n)\|^2+C\tau_c\|U^{n-1}_C-u(t_{n-1})\|^2.
\end{split}
\end{equation}
Replacing $n$ by $j$ and summing from 2 to $n$, we have
\begin{equation}\label{5.7.8}\begin{split}
\mathbb{L}[\xi_c^{n}]-\mathbb{L}[\xi_c^{1}]& \leq 8\tau_c \sum_{j=2}^n(\|R_t^{j-\theta}\|^2+\|\mathcal{D}_{\tau_c}\rho_c^{j-\theta}\|^2+\|E_2^{j-\theta}\|^2+\|E_3^{j-\theta}\|^2)+C\tau_c\sum_{j=1}^n\|U^j_C-u(t_j)\|^2
\\& \leq 8\tau_c \sum_{j=2}^n(\|R_t^{j-\theta}\|^2+\|\mathcal{D}_{\tau_c}\rho_c^{j-\theta}\|^2+\|E_2^{j-\theta}\|^2+\|E_3^{j-\theta}\|^2)+C\tau_c\sum_{j=1}^n(|\xi_c^j|_{\mu}^2+|\rho_c^j|_{\mu}^2).
\end{split}
\end{equation}
By Lemma \ref{5.6}
\begin{equation}\label{5.7.9}\begin{split}
\mathbb{L}[\xi_c^{n}] &\ge  \frac{1}{1-\theta}|\xi_c^n|_{\mu}^2,
\\ \mathbb{L}[\xi_c^{1}]& =  (3-2\theta)|\xi_c^1|^2_{\mu}-(1-2\theta)|\xi_c^0|^2_{\mu}+(2-\theta)(1-2\theta)|\xi_c^1-\xi_c^0|_{\mu}^2 \leq C(|\xi_c^0|^2_{\mu}+|\xi_c^1|^2_{\mu}),
\end{split}\end{equation}
as well as (\ref{5.7.8}) one can derive
\begin{equation}\label{5.7.10}\begin{split}
|\xi_c^{n}|_{\mu}^2& \leq C\tau_c \sum_{j=2}^n(\|R_t^{j-\theta}\|^2+\|\mathcal{D}_{\tau_c}\rho_c^{j-\theta}\|^2+\|E_2^{j-\theta}\|^2+\|E_3^{j-\theta}\|^2)+C\tau_c\sum_{j=1}^n|\rho_c^j|_{\mu}^2
\\ &\quad + C(|\xi_c^0|^2_{\mu}+|\xi_c^1|^2_{\mu})+C\tau_c\sum_{j=1}^n |\xi_c^j|_{\mu}^2.
\end{split}
\end{equation}
Using the discrete Gronwall inequality we have for $n \ge 2$,
\begin{equation}\label{5.7.11}\begin{split}
|\xi_c^{n}|_{\mu}^2& \leq C\tau_c \sum_{j=2}^n(\|R_t^{j-\theta}\|^2+\|E_2^{j-\theta}\|^2+\|E_3^{j-\theta}\|^2)+C \tau_c\sum_{j=2}^n\|\mathcal{D}_{\tau_c}\rho_c^{j-\theta}\|^2+C\tau_c\sum_{j=1}^n|\rho_c^j|_{\mu}^2
\\ &\quad + C(|\xi_c^0|^2_{\mu}+|\xi_c^1|^2_{\mu})
\\ & \leq C(\tau_c n)\tau_c^4 + C \tau_c\sum_{j=2}^n(\|\mathcal{D}_{\tau_c}\rho_c^{j-\theta}-\frac{\partial}{\partial t}\rho_c^{j-\theta}\|^2 + |\frac{\partial}{\partial t}\rho_c^{j-\theta}|_{\mu}^2)+C\tau_c\sum_{j=1}^n|\rho_c^j|_{\mu}^2
\\ &\quad + C(|\xi_c^0|^2_{\mu}+|\xi_c^1|^2_{\mu})
\\ & \leq C t_n \tau_c^4 + C \tau_c\sum_{j=2}^n(\tau_c^4 + h^{2r-2{\mu}}\|u_t^{j-\theta}\|_{H^r}^2)+C\tau_c\sum_{j=1}^n h^{2r-2{\mu}}\|u^{j-\theta}\|_{H^r}^2+ C(|\xi_c^0|^2_{\mu}+|\xi_c^1|^2_{\mu})
\\ & \leq C t_n \tau_c^4 + C t_n \max_{t\in[0,T]}(\|u_t\|_{H^r}^2+\|u\|_{H^r}^2)h^{2r-2{\mu}}+C|U_C^0-u_0|_{\mu}^2+C|\xi_c^1|^2_{\mu}.
\end{split}\end{equation}
To estimate $|\xi_c^1|_{\mu}$, subtract (\ref{5.7.3}) from (\ref{3.9}) to derive
\begin{equation}\label{5.7.12}\begin{split}
(\partial_{\frac{1}{2}}\xi_c,v_h)+B(\xi_c^{\frac{1}{2}},v_h)+(f^{\frac{1}{2}}(U_C)-f^{\frac{1}{2}}(u),v_h)=(E_1-\partial_{\frac{1}{2}} \rho_c -E_2^{\frac{1}{2}}+E_3^{\frac{1}{2}},v_h).
\end{split}
\end{equation}
By taking $v_h=(\xi_c^1-\xi_c^0)/\tau_c$ in (\ref{5.7.12}) and using the formula
\begin{equation}\label{5.7.13}\begin{split}
B(\xi_c^1,\frac{\xi_c^1-\xi_c^0}{\tau_c}) \ge  \frac{1}{2\tau_c}(|\xi_c^1|_{\mu}^2-|\xi_c^0|_{\mu}^2)
\end{split}
\end{equation}
as well as the quite similar analysis, we have
\begin{equation}\label{5.7.14}\begin{split}
|\xi_c^1|_{\mu}^2 \leq C(\tau_c^4+h^{2r-2\mu})+C|U_C^0-u_0|_{\mu}^2.
\end{split}
\end{equation}
Combining (\ref{5.7.11}), (\ref{5.7.14}) with the property of the orthogonal projector $P_h$ and the equivalence of the seminorm $|u|_{\mu}$ and norm $\|u\|_{\mu}$ within $H^{\mu}_0(\Omega)$ due to Lemma \ref{2.6}, we complete the proof of (\ref{5.7.1}).
\\\textbf{Step II:} For (\ref{5.7.2}), we first estimate the error $|u(t_m)-U_I^m|_{\mu}$ on the fine time mesh.
By the notations introduced in (\ref{4.1.15.}), we have
\begin{equation}\label{5.7.15}\begin{split}
U^m_I&=\lambda_m U^{n-1}_C+(1-\lambda_m)U^{n}_C
\\ u(t_m)&=\lambda_m u^{n-1}+(1-\lambda_m)u^{n}+ C \tau_c^2 u_{tt}(\vartheta_m),
\end{split}
\end{equation}
where $\vartheta_m \in (t_{n-1},t_n)$. With (\ref{5.7.15}) and (\ref{5.7.1}) the following result is obvious by triangle inequality
\begin{equation}\label{5.7.16}\begin{split}
|u(t_m)-U_I^m|_{\mu} \leq C(\tau_c^2+h^{r-{\mu}}).
\end{split}
\end{equation}
Next, we replace $n$ by $m$, $\tau_c$ by $\tau$ in (\ref{5.7.4}) respectively and subtract it from (\ref{3.13.}) to get
\begin{equation}\label{5.7.17}\begin{split}
(\mathcal{D}_{\tau}\xi_f^{m-\theta}&,v_h)+B(\xi_f^{m-\theta},v_h)+(1-\theta)(f(U_I^m)+(U_F^m-U_I^m)f_u(U_I^m)-f(u^m),v_h)
\\ &+\theta(f(U_F^{m-1})-f(u^{m-1}),v_h)=(R_t^{m-\theta}-\mathcal{D}_{\tau}\rho_f^{m-\theta}-E^{m-\theta}_2+E^{m-\theta}_3,v_h), \quad \forall v_h \in V_h,
\end{split}
\end{equation}
where $U_F^m-u(t_m) = (U_F^m- P_h u(t_m)) + (P_h u(t_m)- u(t_m)) =\xi_f^m + \rho_f^m$.
\par
Replacing $v_h$ by $\mathcal{D}_{\tau}\xi_f^{m-\theta}$ and recombining the items in (\ref{5.7.17}) to eliminate $(\mathcal{D}_{\tau}\xi_f^{m-\theta},\mathcal{D}_{\tau}\xi_f^{m-\theta})$ by  Cauchy-Schwarz inequality and Young inequality, we have
\begin{equation}\label{5.7.18}\begin{split}
B(\xi_f^{m-\theta}&,\mathcal{D}_{\tau}\xi_f^{m-\theta}) \leq (1-\theta)^2\|f(u^m)-f(U_I^m)-(U_F^m-U_I^m)f_u(U_I^m)\|^2
\\&+ \theta^2 \|f(U_F^{m-1})-f(u^{m-1})\|^2 + \frac{1}{2}\|R_t^{m-\theta}-\mathcal{D}_{\tau}\rho_f^{m-\theta}-E^{m-\theta}_2+E^{m-\theta}_3\|^2.
\end{split}
\end{equation}
Using Taylor expension, we estimate the first item on the rightside of (\ref{5.7.18}) as follows
\begin{equation}\label{5.7.19}\begin{split}
&\|f(u^m)-f(U_I^m)-(U_F^m-U_I^m)f_u(U_I^m)\|
\\&=\|f_u(U_I^m)(u^m-U_I^m)-f_u(U_I^m)(U_F^m-U_I^m)+Cf_{uu}(\eta_m)(u^m-U_I^m)^2\|
\\&=\|f_u(U_I^m)(u^m-U_F^m)+Cf_{uu}(\eta_m)(u^m-U_I^m)^2\|
\\&=\|f_u(U_I^m)(\xi_f^m+\rho_f^m)+Cf_{uu}(\eta_m)(u^m-U_I^m)^2\|
\\&\leq C(\|\xi_f^m\|+\|\rho_f^m\|)+C\|u^m-U_I^m\|_{L^4(\Omega)}^2.
\end{split}
\end{equation}
Combining (\ref{5.7.19}) and (\ref{5.7.18}) with the similar analysis applied to (\ref{5.7.6}), we have for $n \geq 2$
\begin{equation}\label{5.7.20}\begin{split}
|\xi_f^{m}|_{\mu}^2& \leq C  (\tau^4+\tau_c^8+h^{2r-2{\mu}})+C|U_F^0-u_0|_{\mu}^2+C|\xi_f^1|^2_{\mu}.
\end{split}\end{equation}
To estimate $|\xi_f^1|_{\mu}$, subtract (\ref{5.7.3}) from (\ref{3.12.}) to derive
\begin{equation}\label{5.7.21}\begin{split}
(\partial_{\frac{1}{2}}\xi_f&,v_h)+B(\xi_f^{\frac{1}{2}},v_h)+\frac{1}{2}(f(U_I^1)+(U_F^1-U_I^1)f_u(U_I^1)-f(u^1),v_h)
\\ &+\frac{1}{2}(f(U_F^0)-f(u^0),v_h)=(E_1-\partial_{\frac{1}{2}}\rho_f-E^{\frac{1}{2}}_2+E^{\frac{1}{2}}_3,v_h), \quad \forall v_h \in V_h.
\end{split}
\end{equation}
By taking $v_h=(\xi_f^1-\xi_f^0)/\tau$ in (\ref{5.7.21}) and again using the above technique, we can get
\begin{equation}\label{5.7.22}\begin{split}
|\xi_f^1|_{\mu}^2 \leq C(\tau^4+\tau_c^8+h^{2r-2\mu})+C|U_F^0-u_0|_{\mu}^2.
\end{split}
\end{equation}
Combining (\ref{5.7.21}), (\ref{5.7.22}) with the property of the orthogonal projector $P_h$ and the equivalence of the seminorm $|u|_{\mu}$ and norm $\|u\|_{\mu}$ within $H^{\mu}_0(\Omega)$ due to Lemma \ref{2.6}, we complete the proof of (\ref{5.7.2}).
The proof is completed.
\section{Numerical tests}
In this section, we take some numerical examples to test the computational efficiency of TTM method combined with $\theta$-scheme with temporal second-order convergence rate. For implementing the numerical computations in two-dimensional cases, we take rectangular partition for spatial domain $\overline{\Omega}$ and choose continuous bilinear element with basis function $P(x,y)=a+bx+cy+dxy$. In the following numerical tests, we choose three numerical examples based on the space-time domain $[0,1]^2 \times [0,1]$. For convenient implementation of the calculous of the fractional norm $\|\cdot\|_{\mu}$, we use left fractional norm $\|\cdot\|_{J^{\mu}_L(\Omega)}$ instead as they are equivalent within $H_0^{\mu}$. That is to say that we give numerical calculation data by the following norm formula
\begin{equation}\label{6.1.1}\begin{split}
\|u-U_F\|_{J^{\mu}_L(\Omega)}=\Big{(}\|u-U_F\|^2+\|_{RL}D_{a,x}^{\mu} (u-U_F)\|^2+\|_{RL}D_{c,y}^{\mu} (u-U_F)\|^2\Big{)}^{\frac12},
\end{split}
\end{equation}
where
\begin{equation*}\label{6.1.2}\begin{split}
\|_{RL}D_{a,x}^{\mu} (u-U_F)\|^2&=\sum_{e_i}\int_{e_i}\Bigl(_{RL}D_{a,x}^{\mu} u-\frac{1}{\Gamma(1-\mu)}\frac{\partial}{\partial x}\int_a^x \frac{u_{h}(\tau,y) \mathrm{d} \tau}{(x-\tau)^{\mu}}  \Bigr)^2 \mathrm{d} x \mathrm{d} y,
\\\|_{RL}D_{c,y}^{\mu} (u-U_F)\|^2&=\sum_{e_i}\int_{e_i}\Bigl(_{RL}D_{c,y}^{\mu} u-\frac{1}{\Gamma(1-\mu)}\frac{\partial}{\partial y}\int_c^y \frac{u_{h}(x,\tau) \mathrm{d} \tau}{(y-\tau)^{\mu}}  \Bigr)^2 \mathrm{d} x \mathrm{d} y,
\end{split}
\end{equation*}
in which $u_h(x,y)=\textbf{U}_F^T \textbf{N}$ within the element $e$, and $\textbf{N}$ is the element shape function.
\par
At the same time, we also analyze the impact of parameter $M$ on the CPU time and computational accuracy.
\subsection{Numerical data on convergence results}\label{subsec6.1}
\textbf{Example 6.1}
\par
We choose the initial condition $u_0=x^2(x-1)^2y^2(y-1)^2$, the exact solution $u(x,y,t)=e^tx^2(x-1)^2y^2(y-1)^2$, and then the known source term $g(\textbf{z},t)$ can be arrived at.
\par
Table \ref{tab1} mainly shows several cases for different parameters based on the choice of $\tau_{c}=M\tau=\frac{1}{20},M=10$ and $h=1/10,1/20,1/40$: For the fixed $\epsilon=0.01$ and changed fractional parameters $\alpha=1.1,1.5,1.8$, under the case of second-order backward difference time discrete scheme with $\theta=0$, the convergence rate of $\|u-U_F\|$ is approximating the real order $2$ and the convergence orders for errors $\|u-U_F\|_{\mu}$ ($\mu=\frac\alpha2$) are close to the real orders $1.45,1.25,1.10$ (=$2-\mu$), respectively; For the same $\epsilon$ and changed $\alpha=1.3,1.7$, based on the Crank-Nicolson case with $\theta=0.5$, we also get the similar approximation results to the above numerical data with both $\|u-U_F\|$ and $\|u-U_F\|_{\mu}$; Further, by taking $\epsilon=1$, $\theta=0.2,0.4$ and different parameters $\alpha$, we get the same conclusion to the above computation. These numerical results imply that our numerical algorithm is effective. Moreover, for this example, we get the very similar computing accuracy to that calculated by nonlinear Galerkin FE method, which are not provided again because of the same computing results as that listed in Table \ref{tab1}. What's more, from Table \ref{tab1}, one can check that our numerical algorithm can greatly reduce the CPU time.
\par
Table \ref{tab2} continues giving the data statistics on both $\|u-U_F\|$ and $\|u-U_F\|_{\mu}$ with $\tau=\tau_c^2=1/4,1/9,1/16,1/25,1/36$. By choosing $\epsilon=0.1,10$ and different parameters $\theta$ and fractional parameter $\alpha$, we get the same conclusion according to the similar analysis as that discussed for Table \ref{tab1}.
\par
In order to check the temporal convergence order with respect to norm $\|\cdot\|_{\mu}$, we give the calculating results of nonlinear Galerkin FE method with $\epsilon=0.01$ in Table \ref{tab5} by taking $\tau^2=h^{2-\mu}$, which implies that time convergence rate is close to $2$ which is in agreement with the theoretical convergence order of second-order $\theta$ scheme in time.
\par
\begin{figure}[h]
\begin{center}
\begin{minipage}{7.5cm}
  \centering\includegraphics[width=7cm]{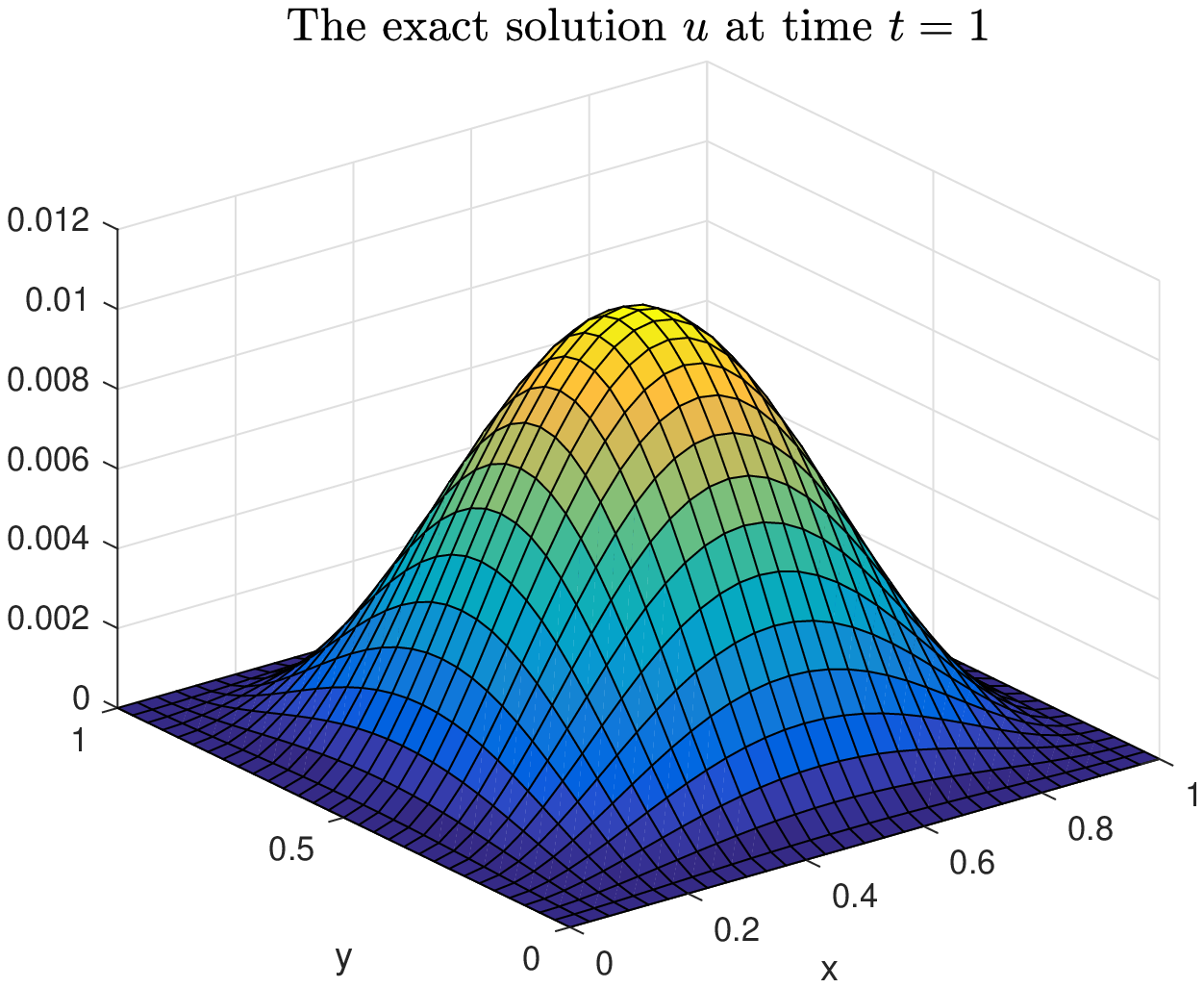}
  \caption{The exact solution $u$ with $\epsilon=1$, $\theta=0.25$, $\alpha=1.4$, $h=\frac{1}{30}$, $\tau=\tau_c/M=\frac{1}{100}$}\label{C1}
\end{minipage}
\begin{minipage}{7.5cm}
  \centering\includegraphics[width=7cm]{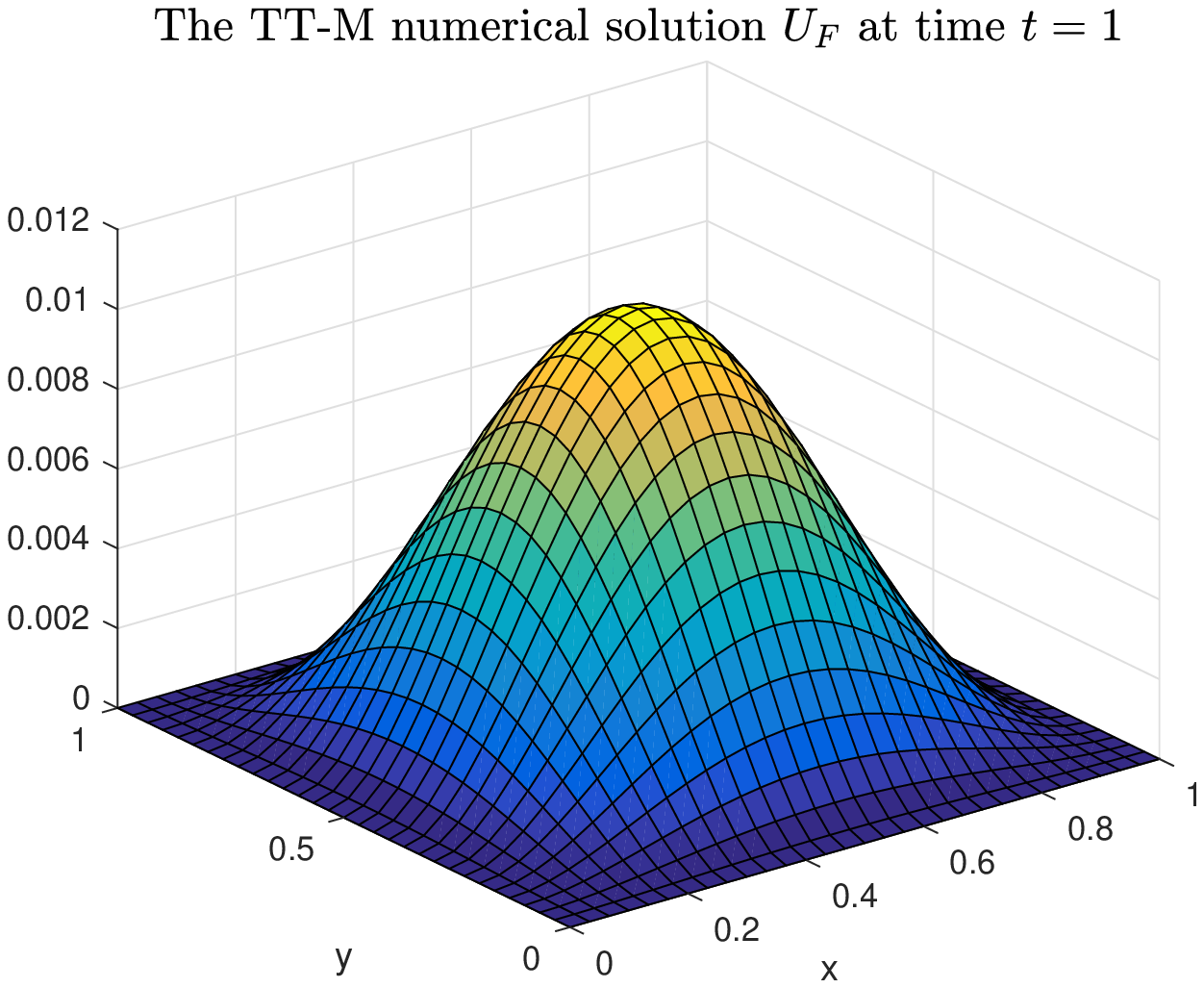}
  \caption{The TT-M FE solution $U_F$ with $\epsilon=1$, $\theta=0.25$, $\alpha=1.4$, $h=\frac{1}{30}$, $\tau=\tau_c/M=\frac{1}{100}$}\label{C2}
\end{minipage}
\end{center}
\end{figure}

To observe the numerical behavior of TT-M FE solution, Figure \ref{C1} and Figure \ref{C2} with $\epsilon=1$, $\theta=0.25$, $\alpha=1.4$, $h=\frac{1}{30}$, $\tau=\tau_c/M=\frac{1}{100}$, $M=10$ show the surfaces for the exact solution $u$ and the TT-M FE solution $U_F$ at time $t=1$, respectively. The result shows that the TT-M FE solution can well approximate the exact solution.
\\
\\
\textbf{Example 6.2}
\par
Now we provide the second example only covering initial condition $u_0=x^2(x-1)^2y^2(y-1)^2$, the source term $g(\textbf{z},t)=0$ and the diffusion coefficient $\epsilon=0.01$. Because of unknown exact solution, we choose the numerical solution with $\tau=h=1/100$ as the approximating exact solution.
Table \ref{tab6} shows the space convergence data in $L^2$-norm with $\tau_c=\sqrt{\tau}=1/10$, changed $h=1/4,1/8,1/16$ and different $\theta$ and $\alpha$. Table \ref{tab7} lists the space-time convergence results in $L^2$-norm containing $h=\tau=\tau_c/2=1/4,1/10,1/20$, which imply the temporal convergence rate is approximating $2$ and is not impacted by the changed parameter $\theta$.
\par
\begin{figure}[h]
\begin{center}
\begin{minipage}{7.5cm}
  \centering\includegraphics[width=7cm]{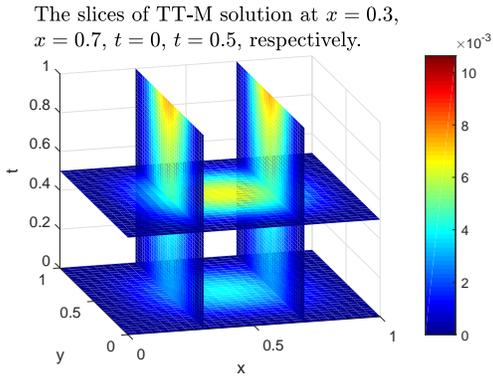}
  \caption{The behavior of TT-M numerical solution with $\theta=0.2$, $\alpha=1.5$, $h=\frac{1}{30}$, $\tau_c=\sqrt{\tau}=\frac{1}{5}$}\label{C3}
\end{minipage}
\begin{minipage}{7.5cm}
  \centering\includegraphics[width=7cm]{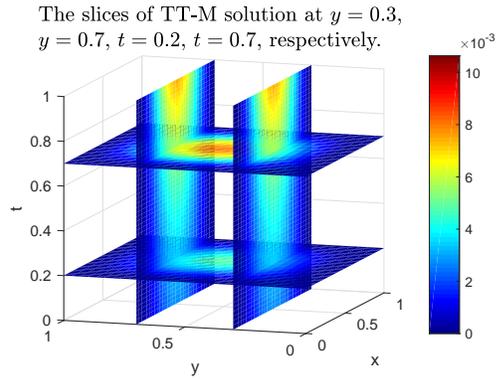}
  \caption{The behavior of TT-M numerical solution with $\theta=0.2$, $\alpha=1.5$, $h=\frac{1}{30}$, $\tau_c=\sqrt{\tau}=\frac{1}{5}$}\label{C4}
\end{minipage}
\end{center}
\end{figure}
\begin{figure}[h]
\begin{center}
\begin{minipage}{7.5cm}
  \centering\includegraphics[width=7cm]{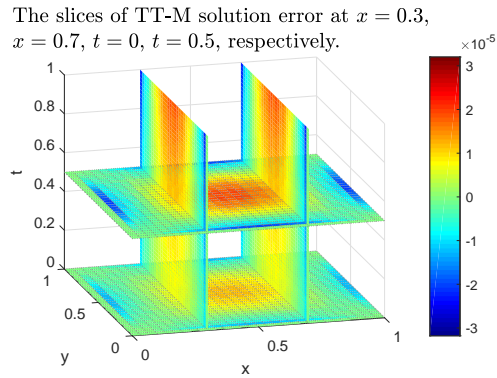}
  \caption{The behavior of error $U_F-u$ with $\theta=0.2$, $\alpha=1.5$, $h=\frac{1}{30}$, $\tau_c=\sqrt{\tau}=\frac{1}{5}$}\label{C5}
\end{minipage}
\begin{minipage}{7.5cm}
  \centering\includegraphics[width=7cm]{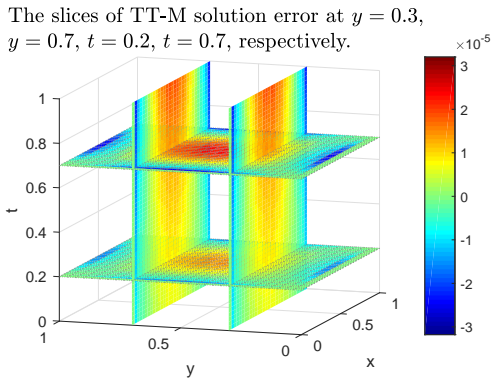}
  \caption{The behavior of error $U_F-u$ with $\theta=0.2$, $\alpha=1.5$, $h=\frac{1}{30}$, $\tau_c=\sqrt{\tau}=\frac{1}{5}$}\label{C6}
\end{minipage}
\end{center}
\end{figure}

For checking the behaviors of numerical solution and error, we consider the numerical performance with $\theta=0.2$, $\alpha=1.5$, $h=\frac{1}{30}$, $\tau_c=\sqrt{\tau}=\frac{1}{5}$ in Figures \ref{C3}-\ref{C6}.
Figure \ref{C3} show the behavior of TT-M numerical solution on different slices at $x=0.3$, $x=0.7$, $t=0$ and $t=0.5$. When taking $x=0.3$ and $x=0.7$, Figure \ref{C3} describes the behavior of TT-M numerical solution on $y,t$ plane, which tell us that the numerical solution may be similar. For the case at $t=0$ and $t=0.5$, the behavior of the TT-M numerical solution on $x,y$ plane shows the value of numerical solution increases with the increase of time. Similarly, Figure
\ref{C4} also shows the numerical behavior of solution based on given slices at $y=0.3$, $y=0.7$, $t=0.2$ and $t=0.7$. From Figures \ref{C3}-\ref{C4}, ones can see the overall trend of the numerical solution based on three space-time parameters $x,y,t$.
\par
For the fixed splice at $x=0.3$ or $x=0.7$, the behavior of error $U_F-u$ in Figure \ref{C5} tells us that the absolute error gradually becomes larger with the increase of time from $t=0$ to $1$ and that how the error relies on the change of variable $y$. We also see the similar behavior of error $U_F-u$ in Figure \ref{C6}.
\\
\\
\\
\textbf{Example 6.3}
\par
\begin{figure}[h]
\begin{center}
\begin{minipage}{7.5cm}
  \centering\includegraphics[width=7cm]{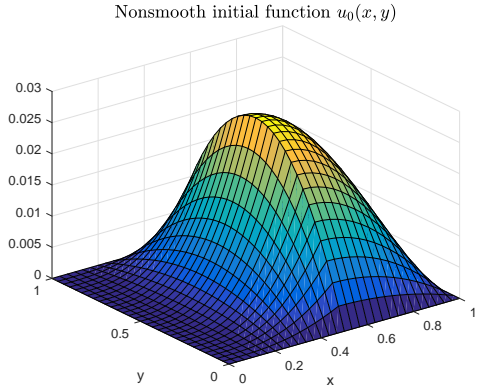}
  \caption{The chosen non-smooth initial function $u_0(x,y)$}\label{C7}
\end{minipage}
\begin{minipage}{7.5cm}
  \centering\includegraphics[width=7cm]{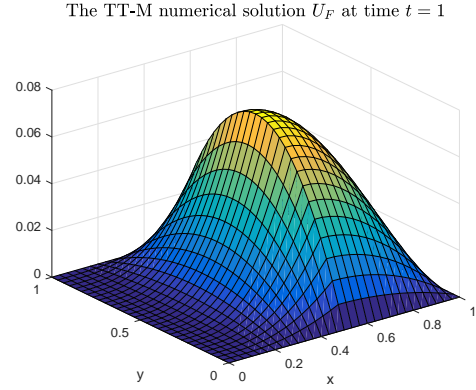}
  \caption{The non-smooth numerical solution with $\theta=0.25$, $\alpha=1.2$, $h=\frac{1}{30}$, $\tau_c=\frac{1}{5}$, $M=4$}\label{C8}
\end{minipage}
\end{center}
\end{figure}

Here, we provide an example with the diffusion coefficient $\epsilon=0.01$, the source term $g(\textbf{z},t)=0$ and non-smooth initial data
\begin{equation}\label{6.1}\begin{split}
u_0(x,y) = \left\{ \begin{array}{rl} x^3(1-x^3)y(1-y), &  x \in [0,0.5],\\ \frac{7}{16}x(1-x)y(1-y), & x \in (0.5,1].
\end{array} \right.
\end{split}
\end{equation}
For this problem, we cannot find the exact solution. So, we need to take the numerical solution under the condition $\tau=\frac{1}{100}$, $h=\frac{1}{100}$ as the approximate exact solution. One can see the detailed numerical data containing errors and convergence rate in Tables \ref{tab8}-\ref{tab9}. By the similar analysis to that in the second example, one can know that our method is also effective for the current example with non-smooth solution. Figure \ref{C7} shows the surface for the given initial data, from which one can easily see that the initial value $u_0(x,y)$ is a non-smooth function. Figure \ref{C8} provides the numerical surface under the condition $\theta=0.25$, $\alpha=1.2$, $h=\frac{1}{30}$, $\tau_c=\frac{1}{5}$, $M=4$, which implies that the numerical solution is non-smooth, and has been impacted by the non-smooth initial data $u_0(x,y)$.
\par
All in all, based on the above three numerical examples with smooth solution, smooth initial data, non-smooth initial function, one can find that our numerical method is effective for solving the nonlinear space fractional Allen-Cahn problem, while saving the computing time (CPU time) with the comparison to that calculated by standard nonlinear Galerkin FE method, and get the second-order time convergence rate which is in agreement with second-order $\theta$ scheme.

\subsection{The impact of $M$ on CPU time and computational accuracy}\label{subsec6.2}
\begin{figure}[h]
\begin{center}
\begin{minipage}{7.5cm}
  \centering\includegraphics[width=7cm]{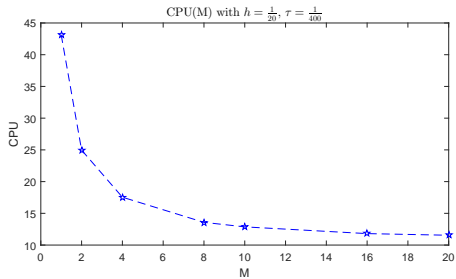}
  \caption{CPU(M) based on $\epsilon=0.1$, $\theta=0.5$, $\alpha=1.2$}\label{C9}
\end{minipage}
\begin{minipage}{7.5cm}
  \centering\includegraphics[width=7cm]{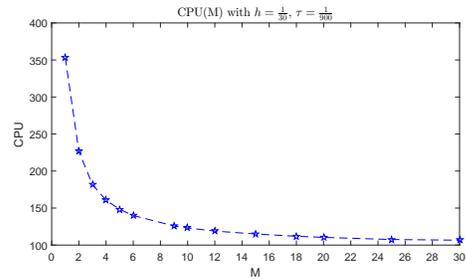}
  \caption{CPU(M) based on $\epsilon=1$, $\theta=0.1$, $\alpha=1.8$}\label{C10}
\end{minipage}
\end{center}
\end{figure}

In subsection \ref{subsec6.1}, we have implemented three numerical examples and given the detailed calculated data analysis for convergence results covering errors, convergence rate and CPU time. For checking the computational efficiency of fast TT-M FE algorithm, we need to consider the impact of parameter $M$ on CPU time. Here, we make the related analysis by choosing only the first example in subsection \ref{subsec6.1}. In Figure \ref{C9}, with the fixed parameters $\epsilon=0.1$, $\theta=0.5$, $\alpha=1.2$, $h=1/20$ and $\tau=1/400$, the distribution point graphs of CPU$(M)$ depending on different $M$ are depicted. From the Figure \ref{C9}, one can clearly see that TT-M FE algorithm needs much less CPU time with $M\geq 2$ than that yielded by nonlinear Galerkin FE methods with $M=1$. Moreover, one can also find that the computing time of fast TT-M FE algorithm gradually reduces when $M$ increases from $2$ to $20$, from which one can know that the most efficient calculation is produced at $M=20$(when $\tau=\tau_c^2=1/400$). At the same time, one can see when $M$ tends to the maximum value $M=1/\tau_c=20$, the CPU(M) changes very slowly. For Figure \ref{C10} with parameters $\alpha=1.8$, $\theta=0.1$ and $\epsilon=1$, $h=1/30$, $\tau=1/900$, we can get the similar conclusion.
\par
\begin{figure}[h]
\begin{center}
\begin{minipage}{7.5cm}
  \centering\includegraphics[width=7cm]{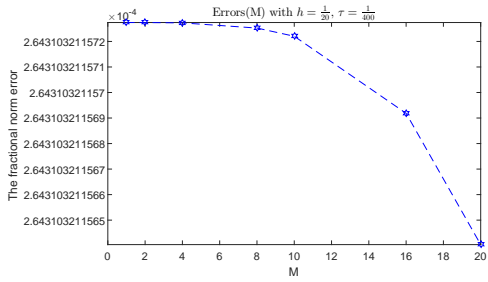}
  \caption{Errors(M) based on $\epsilon=0.1$, $\theta=0.5$, $\alpha=1.2$}\label{C11}
\end{minipage}
\begin{minipage}{7.5cm}
  \centering\includegraphics[width=7cm]{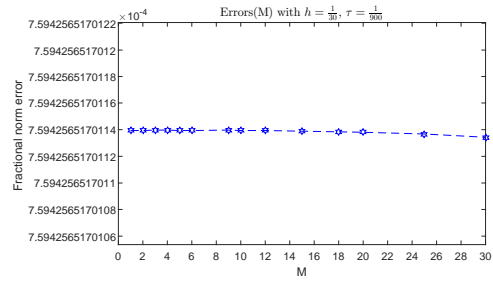}
  \caption{Errors(M) based on $\epsilon=1$, $\theta=0.1$, $\alpha=1.8$}\label{C12}
\end{minipage}
\end{center}
\end{figure}

In what follows, we will check the impact of $M$ on the computational accuracy. With the same fixed parameters chosen as that in Figure \ref{C9}, we compute the error cases Errors(M) based on the changed parameter $M$ in Figure \ref{C11}, from which one clearly see that all errors in norm $\|\cdot\|_{\mu}$ for different parameter $M$ are close to $2.6431032115\times 10^{-4}$. These error data illustrate that the parameter $M$ has very small impact on the computational accuracy. For Figure \ref{C12} with the same parameters as that given in Figure \ref{C10}, we get almost the same results.
\par
 In summary, based on the discussions of the impact of parameter $M$ on both CPU$(M)$ and Errors(M), ones can know that for saving the computing time considerably, a large parameter $M$ (For example $M=\frac{1}{\tau_c}$) may be preferred; whilst any choice of parameter $M$ in $[2,\frac{1}{\tau_c}]$ will arrive at almost the same errors.

\begin{table}[h]
\centering
 \caption{Numerical results of TT-M method with smooth solution and $\tau_c=10\tau=\frac{1}{20}$}\label{tab1}
 \begin{tabular}{ccccccccrr}
  \toprule
    $\epsilon$ & $\theta$ & $\alpha$ & $h$ & $\|u-U_F\|$ &  rate & $\|u-U_F\|_{\mu}$  &  rate  &TT-M(s) &FE(s)\\
\midrule
	&		&	1.1	&	  1/10		&	7.4834E-05	&		    &	6.1736E-04	&	(1.45)  &	1.43 	&	1.14\\
	&		&		&	  1/20		&	1.6390E-05	&	2.191 	&	2.0950E-04	&	1.559 	&	7.87 	&	8.79\\
	&		&		&	  1/40		&	3.6977E-06	&	2.148 	&	7.2897E-05	&	1.523 	&	76.58 	&	139.88\\
	&		&	1.5	&	  1/10		&	7.4835E-05	&		    &	1.3996E-03	&	(1.25)  &	1.29 	&	1.38\\
0.01&	0	&		&	  1/20		&	1.6391E-05	&	2.191 	&	5.6701E-04	&	1.304 	&	7.68 	&	9.14\\
	&		&		&	  1/40		&	3.6986E-06	&	2.148 	&	2.3191E-04	&	1.290 	&	76.38 	&	137.14\\
	&		&	1.8	&	  1/10		&	7.4842E-05	&		    &	2.6543E-03	&	(1.10)  &	1.30 	&	1.39\\
	&		&		&	  1/20		&	1.6403E-05	&	2.190 	&	1.2168E-03	&	1.125 	&	7.77 	&	8.80\\
	&		&		&	  1/40		&	3.7156E-06	&	2.142 	&	5.5813E-04	&	1.124 	&	76.63 	&	137.27\\
\midrule
	&		&	1.3	&	  1/10		&	7.4834E-05	&		    &	9.2218E-04	&	(1.35)  &	1.30 	&	1.38\\
	&		&		&	  1/20		&	1.6390E-05	&	2.191 	&	3.4253E-04	&	1.429 	&	7.82 	&	9.01\\
	&		&		&	  1/40		&	3.6976E-06	&	2.148 	&	1.2940E-04	&	1.404 	&	70.32 	&	140.94\\
0.01&	0.5	&	1.7	&	  1/10		&	7.4837E-05	&		    &	2.1428E-03	&	(1.15)	&	1.34 	&	1.43\\
	&		&		&	  1/20		&	1.6395E-05	&	2.190 	&	9.4353E-04	&	1.183 	&	7.54 	&	9.17\\
	&		&		&	  1/40		&	3.7043E-06	&	2.146 	&	4.1677E-04	&	1.179 	&	70.95 	&	140.56\\
\midrule
	&		&	1.4	&	  1/10		&	8.6946E-05	&		    &	1.0994E-03	&	(1.30)  &	0.84 	&	1.35\\
	&		&		&	  1/20		&	1.9943E-05	&	2.124 	&	4.3079E-04	&	1.352 	&	5.53 	&	8.76\\
	&		&		&	  1/40		&	4.5586E-06	&	2.129 	&	1.7096E-04	&	1.333 	&	71.51 	&	137.86\\
1	&	0.2	&	1.6	&	  1/10		&	9.1050E-05	&		    &	1.6724E-03	&	(1.20)  &	0.83 	&	1.35\\
	&		&		&	  1/20		&	2.1498E-05	&	2.082 	&	7.1513E-04	&	1.226 	&	5.41 	&	8.76\\
	&		&		&	  1/40		&	5.0516E-06	&	2.089 	&	3.0725E-04	&	1.219 	&	71.03 	&	137.86\\
\midrule
	&		&	1.2	&	  1/10		&	8.3180E-05	&		    &	7.3204E-04	&	(1.40)  &	0.77 	&	1.38\\
	&		&		&	  1/20		&	1.8623E-05	&	2.159 	&	2.6203E-04	&	1.482 	&	5.38 	&	8.98\\
	&		&		&	  1/40		&	4.1875E-06	&	2.153 	&	9.5800E-05	&	1.452 	&	72.39 	&	137.75\\
1	&	0.4	&	1.9	&	  1/10		&	1.0233E-04	&		    &	3.1565E-03	&	(1.05)	&	0.79 	&	1.34\\
	&		&		&	  1/20		&	2.5109E-05	&	2.027 	&	1.5322E-03	&	1.043 	&	5.44 	&	8.89\\
	&		&		&	  1/40		&	6.1206E-06	&	2.036 	&	7.3944E-04	&	1.051 	&	71.91 	&	136.88\\

 \bottomrule
 \end{tabular}
\end{table}

\begin{table}[h]
\centering
 \caption{Numerical results of TT-M method with smooth solution and $\tau=h$}\label{tab2}
 \begin{tabular}{cccrrccccrr}
  \toprule
    $\epsilon$ & $\theta$ & $\alpha$ & $\tau_c$ & $\tau$ & $\|u-U_F\|$ &  rate & $\|u-U_F\|_{\mu}$  &  rate  & TT-M(s)& FE(s)\\
\midrule
	&		&		&	  1/2 	&	  1/4 	&	5.2592E-04	&		    &	2.5003E-03	&	(1.45)  &	0.37 	&	0.35\\
	&		&		&	  1/3 	&	  1/9 	&	9.5352E-05	&	2.106 	&	7.2231E-04	&	1.531 	&	0.34 	&	0.62\\
0.1	&	0.1	&	1.1	&	  1/4 	&	  1/16	&	2.7142E-05	&	2.184 	&	2.9317E-04	&	1.567 	&	2.08 	&	2.16\\
	&		&		&	  1/5 	&	  1/25	&	1.0352E-05	&	2.160 	&	1.4740E-04	&	1.541 	&	6.37 	&	7.71\\
	&		&		&	  1/6 	&	  1/36	&	4.7618E-06	&	2.130 	&	8.4775E-05	&	1.517 	&	18.90 	&	25.97\\
\midrule
	&		&		&	  1/2 	&	  1/4 	&	5.2881E-04	&	2.144 	&	4.3232E-03	&	(1.25) 	&	0.04 	&	0.36\\
	&		&		&	  1/3 	&	  1/9 	&	9.7798E-05	&	2.081 	&	1.5691E-03	&	1.250 	&	0.26 	&	0.63\\
0.1	&	0.1	&	1.5	&	  1/4 	&	  1/16	&	2.8645E-05	&	2.134 	&	7.4260E-04	&	1.300 	&	1.08 	&	2.11\\
	&		&		&	  1/5 	&	  1/25	&	1.1130E-05	&	2.118 	&	4.1788E-04	&	1.288 	&	4.79 	&	7.82\\
	&		&		&	  1/6 	&	  1/36	&	5.1645E-06	&	2.106 	&	2.6226E-04	&	1.278 	&	18.65 	&	26.99\\
\midrule
	&		&		&	  1/2 	&	  1/4 	&	5.4173E-04	&	2.118 	&	6.6747E-03	&	(1.10) 	&	0.03 	&	0.36\\
	&		&		&	  1/3 	&	  1/9 	&	1.0644E-04	&	2.007 	&	2.8805E-03	&	1.036 	&	0.23 	&	0.59\\
0.1	&	0.1	&	1.8	&	  1/4 	&	  1/16	&	3.2626E-05	&	2.055 	&	1.5252E-03	&	1.105 	&	1.06 	&	2.14\\
	&		&		&	  1/5 	&	  1/25	&	1.3002E-05	&	2.062 	&	9.2989E-04	&	1.109 	&	4.73 	&	7.44\\
	&		&		&	  1/6 	&	  1/36	&	6.1325E-06	&	2.061 	&	6.2092E-04	&	1.108 	&	18.69 	&	27.37\\
\midrule
	&		&		&	  1/2 	&	  1/4 	&	5.2540E-04	&	2.026 	&	3.2647E-03	&	(1.35) 	&	0.03 	&	0.35\\
	&		&		&	  1/3 	&	  1/9 	&	9.5515E-05	&	2.102 	&	1.0570E-03	&	1.391 	&	0.26 	&	0.62\\
0.1	&	0.3	&	1.3	&	  1/4 	&	  1/16	&	2.7370E-05	&	2.172 	&	4.6337E-04	&	1.433 	&	1.07 	&	2.14\\
	&		&		&	  1/5 	&	  1/25	&	1.0501E-05	&	2.146 	&	2.4662E-04	&	1.413 	&	4.75 	&	7.36\\
	&		&		&	  1/6 	&	  1/36	&	4.8451E-06	&	2.121 	&	1.4826E-04	&	1.395 	&	18.89 	&	26.02\\
\midrule
	&		&		&	  1/2 	&	  1/4 	&	5.3454E-04	&	2.141 	&	5.7670E-03	&	(1.15) 	&	0.02 	&	0.34\\
	&		&		&	  1/3 	&	  1/9 	&	1.0223E-04	&	2.040 	&	2.3497E-03	&	1.107 	&	0.25 	&	0.62\\
0.1	&	0.3	&	1.7	&	  1/4 	&	  1/16	&	3.0772E-05	&	2.087 	&	1.1993E-03	&	1.169 	&	1.09 	&	2.12\\
	&		&		&	  1/5 	&	  1/25	&	1.2121E-05	&	2.088 	&	7.1226E-04	&	1.168 	&	4.71 	&	7.38\\
	&		&		&	  1/6 	&	  1/36	&	5.6694E-06	&	2.084 	&	4.6600E-04	&	1.164 	&	18.73 	&	27.51\\
\midrule
	&		&		&	  1/2 	&	  1/4 	&	5.5713E-04	&	2.088 	&	3.6738E-03	&	(1.30) 	&	0.03 	&	0.35\\
	&		&		&	  1/3 	&	  1/9 	&	1.0855E-04	&	2.017 	&	1.2684E-03	&	1.311 	&	0.23 	&	0.60\\
10	&	0.2	&	1.4	&	  1/4 	&	  1/16	&	3.2081E-05	&	2.118 	&	5.8170E-04	&	1.355 	&	1.05 	&	2.19\\
	&		&		&	  1/5 	&	  1/25	&	1.2395E-05	&	2.131 	&	3.1948E-04	&	1.343 	&	4.73 	&	7.46\\
	&		&		&	  1/6 	&	  1/36	&	5.7025E-06	&	2.129 	&	1.9659E-04	&	1.332 	&	18.81 	&	27.09\\
\midrule
	&		&		&	  1/2 	&	  1/4 	&	5.5390E-04	&	2.083 	&	7.7054E-03	&	(1.05) 	&	0.03 	&	0.34\\
	&		&		&	  1/3 	&	  1/9 	&	1.4809E-04	&	1.627 	&	3.5389E-03	&	0.960 	&	0.25 	&	0.61\\
10	&	0.5	&	1.9	&	  1/4 	&	  1/16	&	4.4993E-05	&	2.071 	&	1.9440E-03	&	1.041 	&	1.07 	&	2.14\\
	&		&		&	  1/5 	&	  1/25	&	1.9066E-05	&	1.924 	&	1.2154E-03	&	1.052 	&	4.72 	&	7.41\\
	&		&		&	  1/6 	&	  1/36	&	8.7374E-06	&	2.140 	&	8.2780E-04	&	1.053 	&	18.71 	&	27.22\\

 \bottomrule
 \end{tabular}
\end{table}

\begin{table}[h]
\centering
 \caption{Convergence results of FM method covering smooth solution and $\tau^2=h^{2-\mu}$}\label{tab5}
 \begin{tabular}{ccrrcc}
  \toprule
     $\theta$ & $\alpha$ & $h$ & $\tau$ & $\|u-U_F\|_{\mu}$  &  rate\\
\midrule
		&		&	  1/3 	&	  1/2 	&	4.9870E-03	&		\\
0.25	&	1.4	&	  1/12	&	  1/5 	&	8.5823E-04	&	1.920 	\\
		&		&	  1/40	&	  1/11	&	1.7100E-04	&	2.046 	\\

 \bottomrule
 \end{tabular}
\end{table}
\begin{table}[h]
\centering
 \caption{TT-M convergence results including smooth initial value and $\epsilon=0.01$}\label{tab6}
 \begin{tabular}{ccccrcc}
  \toprule
 $\theta$ & $\alpha$ & $\tau_c$ & $\tau$ & $h$& $\|u-U_F\|$  &  rate\\
\midrule
    0.2	&	1.5	&	 1/10	    &	   1/100	&	  1/4 	&	5.2001E-04	&		\\
		&		&				&		        &	  1/8 	&	1.2057E-04	&	2.109 	\\
		&		&				&		        &	  1/16	&	2.5975E-05	&	2.215 	\\
	0.5	&	1.8	&	 1/10	    &	   1/100	&	  1/4 	&	5.1746E-04	&		\\
		&		&				&		        &	  1/8 	&	1.1877E-04	&	2.123 	\\
		&		&				&		        &	  1/16	&	2.4981E-05	&	2.249 	\\
	0	&	1.2	&	 1/10	    &	   1/100	&	  1/4 	&	5.2135E-04	&		\\
		&		&				&		        &	  1/8 	&	1.2140E-04	&	2.102 	\\
		&		&				&		        &	  1/16	&	2.6468E-05	&	2.197 	\\
 \bottomrule
 \end{tabular}
\end{table}

\begin{table}[h]
\centering
 \caption{TT-M convergence results including smooth initial value}\label{tab7}
 \begin{tabular}{ccrrrcc}
  \toprule
 $\theta$ & $\alpha$& $\tau_c$& $\tau$ & $h$  & $\|u-U_F\|$  &  rate\\
\midrule
	0.2	&	1.5	&	  1/2 	&	  1/4 	&	  1/4 	&	5.2208E-04	&		\\
		&		&	  1/5 	&	  1/10	&	  1/10	&	7.4467E-05	&	2.125 	\\
		&		&	  1/10	&	  1/20	&	  1/20	&	1.5748E-05	&	2.241 	\\
	0.5	&	1.8	&	  1/2 	&	  1/4 	&	  1/4 	&	5.1792E-04	&		\\
		&		&	  1/5 	&	  1/10	&	  1/10	&	7.2547E-05	&	2.145 	\\
		&		&	  1/10	&	  1/20	&	  1/20	&	1.4837E-05	&	2.290 	\\
	0	&	1.2	&	  1/2 	&	  1/4 	&	  1/4 	&	5.2478E-04	&		\\
		&		&	  1/5 	&	  1/10	&	  1/10	&	7.5598E-05	&	2.115 	\\
		&		&	  1/10	&	  1/20	&	  1/20	&	1.6253E-05	&	2.218 	\\

 \bottomrule
 \end{tabular}
\end{table}

\begin{table}[h]
\centering
 \caption{TT-M numerical results with non-smooth initial data}\label{tab8}
 \begin{tabular}{ccccrcc}
  \toprule
 $\theta$ & $\alpha$ & $\tau_c$ & $\tau$ & $h$& $\|u-U_F\|$  &  rate\\
\midrule
0.25	&	1.2	&	   1/10 	&	   1/100	&	  1/4 	&		2.1641E-03	&		\\
	&		&		&		&	  1/8 	&		4.5377E-04	&	2.254 	\\
	&		&		&		&	  1/16	&		1.0485E-04	&	2.114 	\\
0	&	1.8	&	   1/10 	&	  1/100	&	  1/4 	&		2.1566E-03	&		\\
	&		&		&		&	  1/8 	&		4.5850E-04	&	2.234 	\\
	&		&		&		&	  1/16	&		1.1983E-04	&	1.936 	\\
 \bottomrule
 \end{tabular}
\end{table}

\begin{table}[h]
\centering
 \caption{TT-M numerical results with non-smooth initial data}\label{tab9}
 \begin{tabular}{ccrrrcc}
  \toprule
 $\theta$ & $\alpha$& $\tau_c$& $\tau$ & $h$  & $\|u-U_F\|$  &  rate\\
\midrule
0.25	&	1.2	&	  1/2 	&	  1/4 	&	  1/4 	&		2.1845E-03	&		\\
	&		&	  1/5 	&	  1/10	&	  1/10	&		2.8803E-04	&	2.211 	\\
	&		&	  1/10	&	  1/20	&	  1/20	&		6.6713E-05	&	2.110 	\\
0	&	1.8	&	  1/2 	&	  1/4 	&	  1/4 	&		2.1939E-03	&		\\
	&		&	  1/5 	&	  1/10	&	  1/10	&		3.0113E-04	&	2.167 	\\
	&		&	  1/10	&	  1/20	&	  1/20	&		8.3815E-05	&	1.845 	\\
 \bottomrule
 \end{tabular}
\end{table}
\section{Conclusion}
In this paper, we apply the fast TT-M FE algorithm to nonlinear space fractional Allen-Cahn equations. It is the first time that TT-M FE algorithm is combined with second-order $\theta$ scheme to formulate the fast computing scheme with the detailed analysis of both stability and error estimates for nonlinear space fractional problems. Finally, we choose three examples with smooth solution, smooth initial value and non-smooth initial data to verify our theoretical results, provide the analysis of comparison of CPU time, and give the discussions for the impact of parameter $M$.
In our other works, as talked in the section of Conclusions and future advancements in \cite{Liu1}, the new space-time two-mesh (S-TT-M) method, which is formulated by combining fast TT-M algorithm with Xu's space two-grid method, can be applied to solving nonlinear evolution equations.

\section*{Acknowledgements} This work is supported by the National Natural Science Fund (11661058, 11761053, 11501311), Natural Science Fund of Inner Mongolia Autonomous Region (2016MS0102, 2017MS0107), program for Young Talents of Science and Technology in Universities of Inner Mongolia
Autonomous Region (NJYT-17-A07).
\section*{Conflict of Interest:} The authors declare that they have no conflict of interest.

\end{document}